\documentclass[12pt]{article}
\usepackage{amssymb,amsmath}
\renewcommand{\P}{{\mathcal P}}
\newcommand{\Z}{{\mathbb Z}}
\newcommand{\T}{{\mathcal T}}
\newcommand{\A}{{\mathcal A}}
\newcommand{\G}{{\mathcal G}}
\newcommand{\z}{{\mathcal Z}}
\renewcommand{\H}{{\mathcal H}}

\newcommand{\Q}{{\mathbb Q}}

\renewcommand{\O}{{\cal O}}
\newcommand{\N}{{\rm N}}
\newcommand{\length}{{\rm length}}
\newcommand{\ext}{{\rm ext}}
\newcommand{\id}{{\rm id}}
\newcommand{\Gal}{{\rm Gal}}
\newcommand{\Aut}{{\rm Aut}}

\newcommand{\Tr}{{\rm Tr}}
\newcommand{\ubar}{\overline{u}}
\newcommand{\vbar}{\overline{v}}
\newcommand{\jhat}{\hat{\jmath}}

\newcommand{\pibar}{\overline{\pi}}

\newcommand{\zetabar}{\overline{\zeta}}
\newcommand{\proof}{\noindent{\em Proof: }}
\newcommand{\qed}{\hspace{\fill}$\square$}
\newcommand{\ra}{\rightarrow}
\newcommand{\lra}{\longrightarrow}
\newcommand{\dst}{\displaystyle}

\newcommand{\lbs}[1]
{\begin{picture}(3,13)
\put (11,-3){\line(-1,1){19}}
\put (-5,-1){$\scriptstyle#1$}
\end{picture}}
\newcommand{\ubs}[1]
{\begin{picture}(3,13)
\put (11,-3){\line(-1,1){19}}
\put (3,8){$\scriptstyle#1$}
\end{picture}}
\newcommand{\ufs}[1]
{\begin{picture}(3,13)
\put (-8,-3){\line(1,1){19}}
\put (-4,8){$\scriptstyle#1$}
\end{picture}}
\newcommand{\lfs}[1]
{\begin{picture}(3,13)
\put (-8,-3){\line(1,1){19}}
\put (3,1){$\scriptstyle#1$}
\end{picture}}
\newcommand{\dashline}
{\begin{picture}(2,10)
\put (1,-5){\line(0,1){4}}
\put (1,3){\line(0,1){4}}
\put (1,11){\line(0,1){4}}
\end{picture}}

\newtheorem{theorem}{Theorem}
\newtheorem{lemma}[theorem]{Lemma}
\newtheorem{prop}[theorem]{Proposition}
\newtheorem{cor}[theorem]{Corollary}

\numberwithin{equation}{section}
\numberwithin{theorem}{section}

\title{Extensions of local fields and truncated power series}
\author{Kevin Keating \\
Department of Mathematics \\
University of Florida \\
Gainesville, FL 32611 \\
USA \\[.2cm]
{\tt keating@math.ufl.edu}}
\date{}
\begin{document}

\maketitle

\begin{abstract}
\noindent
Let $K$ be a finite tamely ramified extension of $\Q_p$ and
let $L/K$ be a totally ramified $(\Z/p^n\Z)$-extension.  Let
$\pi_L$ be a uniformizer for $L$, let $\sigma$ be a
generator for $\Gal(L/K)$, and let $f(X)$ be an element of
$\O_K[X]$ such that $\sigma(\pi_L)=f(\pi_L)$.  We show that
the reduction of $f(X)$ modulo the maximal ideal of $\O_K$
determines a certain subextension of $L/K$ up to isomorphism.
We use this result to study the field extensions generated
by periodic points of a $p$-adic dynamical system.
\end{abstract}

\section{Introduction}

     Let $p$ be a prime and let $\Q_p$ denote the $p$-adic
numbers.  In what follows all extensions of $\Q_p$ are
contained in a fixed algebraic closure $\Q_p^{alg}$ of $\Q_p$.
Let $K$ be a finite extension of $\Q_p$ with ramification
index $e$, let $\O_K$ denote the
ring of integers of $K$, and let $\P_K$ denote the maximal
ideal of $\O_K$.  Let $L/K$ be a totally ramified cyclic
extension of degree $p^n$.  Then the residue field
$k=\O_K/\P_K$ of $K$ may be identified with a
subring of $\O_L/\P_L^{ep^n}$ using the Teichm\"uller lifting.
Let $\sigma$ be a generator for $\Gal(L/K)$ and let $\pi_L$
be a uniformizer for $L$.  Then there is a unique
$h_{\pi_L}^{\sigma}(X)\in k[X]/(X^{ep^n})$
such that $\sigma(\pi_L)\equiv\pi_Lh_{\pi_L}^{\sigma}(\pi_L)
\pmod{\P_L^{ep^n+1}}$.  The aim of this paper is to prove the
following:

\begin{theorem} \label{main}
Let $p>3$ and let $K$ be a finite tamely ramified extension
of $\Q_p$ with ramification index $e$.  Let $L/K$ and $L'/K$
be totally ramified $(\Z/p^n\Z)$-extensions such that $L/K$
is contained in a $\Z_p$-extension $L_{\infty}/K$.  Assume:
\begin{quote} (*)
$\left\{\parbox{12.5cm}{There are
generators $\sigma$, $\sigma'$ for $\Gal(L/K)$, $\Gal(L'/K)$
and uniformizers $\pi_L$, $\pi_{L'}$ for $L$, $L'$
such that $h_{\pi_L}^{\sigma}=h_{\pi_{L'}}^{\sigma'}$.}\right.$
\end{quote}
Let $m_0$ be the largest integer such that
$\psi_{L/K}((m_0+1+\frac{1}{p-1})e)<ep^n$.
Then there is $\omega\in\Gal(\Q_p^{alg}/\Q_p)$ such
that $\omega(K)=K$, $\omega$ induces the identity on $k$,
and ${[L\cap\omega(L'):K]}\geq p^{m_0}$.
\end{theorem}

     The function $\psi_{L/K}:[-1,\infty)\ra[-1,\infty)$ and
its inverse $\phi_{L/K}$ are the Hasse-Herbrand functions of
higher ramification theory.  The basic properties of these
functions can be found in Chapters IV and V of \cite{cl} for
finite Galois extensions, and in the appendix of \cite{D} for
finite separable extensions.  We will make frequent use of
the formulas $\psi_{M/K}=\psi_{M/L}\circ\psi_{L/K}$ and
$\phi_{M/K}=\phi_{L/K}\circ\phi_{M/L}$ for finite separable
extensions $M\supset L\supset K$.

     If $K$ contains no primitive $p$th roots of unity then
it can be shown using class field theory that $L/K$ is
contained in a $\Z_p$-extension (see Lemma~\ref{Zp}).  In any
case, Theorem~\ref{main} is valid if either $L/K$ or
$L'/K$ is contained in a $\Z_p$-extension.  If neither 
of $L/K$, $L'/K$ is contained in a $\Z_p$-extension, we still
have the following result:

\begin{theorem} \label{proot}
Let $p>3$ and let $K$ be a finite tamely ramified extension
of $\Q_p$ with ramification index $e$.  Let $L/K$ and $L'/K$
be totally ramified $\Z/p^n\Z$-extensions which satisfy (*).
Then there is $\omega\in\Gal(\Q_p^{alg}/\Q_p)$ such
that $\omega(K)=K$, $\omega$ induces the identity on $k$, and
$[L\cap\omega(L'):K]\geq p^{m_0-1}$.
\end{theorem}

     Suppose $p>3$ and $K/\Q_p$ is unramified.  Then 
$m_0=n-1$ and $K$ contains no primitive $p$th roots of unity.
Furthermore, any automorphism of $\Q_p^{alg}$
which induces the identity on $k$ also induces the identity
on $K$ and hence maps $L$ onto itself.  Therefore we get a simpler
version of Theorem~\ref{main} in this case.

\begin{cor} \label{special}
Let $p>3$, let $K$ be a finite unramified extension of
$\Q_p$, and let $L/K$, $L'/K$ be totally ramified
$(\Z/p^n\Z)$-extensions which satisfy (*).  Then
$[L\cap L':K]\geq p^{n-1}$.
\end{cor}

     Our proof of Theorem~\ref{main} is motivated by
Wintenberger's proof of \cite[Th.\:2]{WZp}, but uses
Deligne's theory of extensions of truncated valuation rings
in place of the field of norms.  In Section~\ref{norm} we present a
slightly modified version of Wintenberger's theorem and use
it to prove a result which is related to Theorem~\ref{main}.
In Section~\ref{tlr} we give an outline of the theory of
truncated local rings based on \cite{D}.  In
Section~\ref{appear} we prove a version of Theorem~\ref{main}
for cyclotomic extensions.  In Section~\ref{proof} we use
this special case to prove the theorem in general, and
show how the same methods can be used to prove
Theorem~\ref{proot}.  In Section~\ref{dynam} we use a variant
of Theorem~\ref{main} to study the field extensions generated
by periodic points of a $p$-adic dynamical system.

\section{The field of norms} \label{norm}

     In \cite{Wab} and \cite{WZp} Wintenberger describes a
remarkable correspondence between groups of power series over
fields of characteristic $p$ and $\Z_p$-extensions of local
fields.  Theorem~\ref{main}
may be viewed as a finite-level version of a part of
this correspondence.  In this section we describe the
connection between Wintenberger's results and Theorem~\ref{main}.

     We begin by recalling the construction of the field of
norms in a special case \cite{cn}.  We define a local field
to be a field complete with respect to a discrete valuation
which has finite residue field.  Let $L_0$ be a local
field whose residue field $k$ has characteristic $p$
and let $L_{\infty}/L_0$ be a totally ramified
$\Z_p$-extension.  For $n\ge0$ let $L_n/L_0$ denote the
subextension of $L_{\infty}/L_0$ of degree $p^n$, let $\O_n$
denote the ring of integers of $L_n$, and let $\P_n$ denote
the maximal ideal of $\O_n$.  Set $r_n=\lceil(p-1)i_n/p\rceil$,
where $i_n$ is the unique (upper and lower) ramification break
of the $(\Z/p\Z)$-extension $L_{n+1}/L_n$.  It follows from
\cite[Prop.\:2.2.1]{cn} that the norm $\N_{L_{n+1}/L_n}$
induces a ring homomorphism $\overline{\N}_{n+1,n}$ from
$\O_{n+1}/(\P_{n+1}^{r_{n+1}})$ onto $\O_n/(\P_n^{r_n})$.
The ring $A_{L_0}(L_{\infty})$ is defined to be the inverse
limit of the rings $\O_n/(\P_n^{r_n})$ with respect to the
maps $\overline{\N}_{n+1,n}$.
Since $\O_n/(\P_n^{r_n})\cong k[X]/(X^{r_n})$ and
$\lim_{n\ra\infty}r_n=\infty$ we have
$A_{L_0}(L_{\infty})\cong k[[X]]$.  The field of norms
$X_{L_0}(L_{\infty})$ of the extension $L_{\infty}/L_0$ is
defined to be the field of fractions of $A_{L_0}(L_{\infty})$.

     We define a compatible sequence of uniformizers for
$L_{\infty}/L_0$ to be a sequence $(\pi_n)_{n\ge0}$ such that
$\pi_n$ is a uniformizer for $L_n$ and
$\N_{L_{n+1}/L_n}(\pi_{n+1})=\pi_n$ for every $n\ge0$.
Associated to each compatible system of uniformizers for
$L_{\infty}/L_0$ we get a uniformizer $(\pibar_n)_{n\ge0}$
for $X_{L_0}(L_{\infty})$, where $\pibar_n$ denotes the image
of $\pi_n$ in $\O_n/\P_n^{r_n}$.  By \cite[Prop.\:2.3.1]{cn}
this construction gives a bijection between the set of
compatible sequences of uniformizers for $L_{\infty}/L_0$ and
the set of uniformizers for $X_{L_0}(L_{\infty})$.

     Let $\sigma\in\Gal(L_{\infty}/L_0)$.  Then for each
$n\ge0$ there is a unique $g_n(X)\in k[X]$ of degree~$<r_n$
such that
\begin{equation} \label{pig}
\frac{\sigma\pi_n}{\pi_n}\equiv g_n(\pi_n)\pmod{\pi_n^{r_n}},
\end{equation}
where we identify $k$ with a subring of $\O_n/\P_n^{r_n}$
using the Teichm\"uller lifting.  If $n\ge1$ we may apply
$\overline{\N}_{n,n-1}$ to (\ref{pig}).  Since
$\overline{\N}_{n,n-1}$ is a ring homomorphism and
$\Gal(L_n/L_0)$ is abelian we get
\begin{equation}
\frac{\sigma\pi_{n-1}}{\pi_{n-1}}\equiv g_n^{\phi}(\pi_{n-1})
\pmod{\pi_{n-1}^{r_{n-1}}},
\end{equation}
where $g_n^{\phi}(X)$ denotes the image of $g_n(X)$ under the
automorphism of $k[X]$ induced by the $p$-Frobenius of $k$.
It follows that
\begin{equation} \label{normcong}
g_{n-1}(X)\equiv g_n^{\phi}(X)\pmod{X^{r_{n-1}}}.  
\end{equation}
Therefore there is $g_{\sigma}(X)\in k[[X]]$ such that
\begin{equation} \label{gsig}
\frac{\sigma\pi_n}{\pi_n}\equiv g_{\sigma}^{\phi^{-n}}(\pi_n)
\pmod{\pi_n^{r_n}}
\end{equation}
for all $n\ge0$.  We define a $k$-linear action of 
$\Gal(L_{\infty}/L_0)$ on $X_{L_0}(L_{\infty})\cong k(\!(X)\!)$
by setting $\sigma\cdot X=Xg_{\sigma}(X)$.

     Let $\A(k)$ denote the set of power series in $k[[X]]$
whose leading term has degree 1.  Then $\A(k)$ with the
operation of substitution forms a group.  The map which
carries $\sigma\in\Aut_k(k(\!(X)\!))$ onto $\sigma(X)\in\A(k)$
gives an isomorphism between $\Aut_k(k(\!(X)\!))$ and
$\A(k)^{op}$.  The subgroup
${\mathcal N}(k)$ of $\A(k)$ consisting of power series with
leading term $X$ is a pro-$p$ group known as the Nottingham group
\cite{cam}.  Let $\Gamma_{L_{\infty}/L_0}^{(\pi_n)}$ denote
the subgroup of $\A(k)$ consisting of power series of the
form $Xg_{\sigma}(X)$ that arise from
elements $\sigma\in\Gal(L_{\infty}/L_0)$ using the compatible
sequence of uniformizers $(\pi_n)$ for $L_{\infty}/L_0$.  Then
$\Gamma_{L_{\infty}/L_0}^{(\pi_n)}$ is isomorphic to $\Z_p$.
The subgroup
$\Gamma_{L_{\infty}/L_0}^{(\pi_n)}$ of $\A(k)$ is determined
up to conjugation by $L_{\infty}/L_0$, and any subgroup of
$\A(k)$ which is conjugate to $\Gamma_{L_{\infty}/L_0}^{(\pi_n)}$
is equal to $\Gamma_{L_{\infty}/L_0}^{(\tilde{\pi}_n)}$ for some 
compatible sequence of uniformizers $(\tilde{\pi}_n)$ for
$L_{\infty}/L_0$.

     Let $K$, $K'$ be local fields with residue field $k$
and let $L/K$, $L'/K'$ be totally ramified extensions.
We say that $L/K$ is $k$-isomorphic to $L'/K'$ if there is
an isomorphism $\tau:L \ra L'$ such that $\tau(K)=K'$ and
$\tau$ induces the identity on $k$.  In this case we write
$L/K\cong_k L'/K'$.  Let $\z(k)$ denote the set of
$k$-isomorphism classes of totally ramified $\Z_p$-extensions
$L_{\infty}/L_0$ such that $L_0$ is a local field with
residue field
$k$.  We put a metric on $\z(k)$ by defining the distance
between the classes $[L_{\infty}/L_0]$ and $[L_{\infty}'/L_0']$
to be $2^{-m}$, where $0\le m\le\infty$ is the largest value
such that $L_m/L_0\cong_k L_m'/L_0'$, and
$m=-1$ if $L_0$ is not $k$-isomorphic to $L_0'$.
Let $\G(k)$ denote the set of conjugacy classes $[\Gamma]$
of subgroups of $\A(k)$ which are isomorphic to $\Z_p$.
We put a metric on $\G(k)$ by defining the distance between
$[\Gamma]$ and $[\Gamma']$ to be $2^{-m}$, where $m$ is
the largest integer such that
$h\Gamma h^{-1}\equiv\Gamma'\pmod{X^{m+1}}$ for some $h\in\A(k)$.

     Since $\Gamma_{L_{\infty}/L_0}^{(\pi_n)}$ is determined
up to conjugacy by the $k$-isomorphism class of
$L_{\infty}/L_0$, we denote its conjugacy
class by $[\Gamma_{L_{\infty}/L_0}]$.  The following result
is essentially a special case of \cite[Cor.\:1.3]{lms}.

\begin{prop} \label{bij}
The map $\Phi:\z(k)\ra\G(k)$ defined by
$\Phi([L_{\infty}/L_0])=[\Gamma_{L_{\infty}/L_0}]$ is a
continuous bijection.
\end{prop}

\proof Since $\lim_{n\ra\infty}r_n=\infty$, the map $\Phi$
is continuous.  To show that $\Phi$ is onto choose
$[\Gamma]\in\G(k)$.  By \cite[Th.\:1]{WZp}
there is a totally ramified $\Z_p$-extension
$L_{\infty}/L_0$ of local fields with residue field $k$ and a
field isomorphism $f:k(\!(X)\!)\ra X_{L_0}(L_{\infty})$ which
induces an isomorphism between $\Gamma$ and the subgroup of
$\Aut(X_{L_0}(L_{\infty}))$ induced by
$\Gal(L_{\infty}/L_0)$.  We write
\begin{equation}
f:(k(\!(X)\!),\Gamma)\overset{\sim}{\lra}(X_{L_0}(L_{\infty}),\Gal(L_{\infty}/L_0)).
\end{equation}
Let $\omega$ be the automorphism of $k$
induced by $f$ and let $\Omega$ be an automorphism of the
separable closure of $L_{\infty}$ which induces $\omega$ on
$k$.  Let $L_0'=\Omega^{-1}(L_0)$ and
$L_{\infty}'=\Omega^{-1}(L_{\infty})$.  Then $\Omega^{-1}$
induces an isomorphism $\Upsilon:X_{L_0}(L_{\infty})\ra
X_{L_0'}(L_{\infty}')$, and
\begin{equation}
\Upsilon\circ f:(k(\!(X)\!),\Gamma)\lra
(X_{L_0'}(L_{\infty}'),\Gal(L_{\infty}'/L_0'))
\end{equation}
is a $k$-linear isomorphism.  It follows that
$\Phi([L_{\infty}'/L_0'])=[\Gamma]$.

     To show that $\Phi$ is one-to-one suppose
$[\Gamma_{L_{\infty}/L_0}]=[\Gamma_{L_{\infty}'/L_0'}]$.
Then there are compatible sequences of uniformizers $(\pi_n)$
for $L_{\infty}/L_0$ and $(\pi_n')$ for $L_{\infty}'/L_0'$
such that $\Gamma_{L_{\infty}/L_0}^{(\pi_n)}=
\Gamma_{L_{\infty}'/L_0'}^{(\pi_n')}$.  Therefore there
is an isomorphism
\begin{equation}
f:(X_{L_0}(L_{\infty}),\Gal(L_{\infty}/L_0))\lra
(X_{L_0'}(L_{\infty}'),\Gal(L_{\infty}'/L_0'))
\end{equation}
which maps $(\pibar_n)$ to $(\pibar_n')$ and induces the identity
on $k$.  It follows from \cite[Th.\:2]{WZp} that $f$
is induced by a $k$-isomorphism from $L_{\infty}/L_0$ to
$L_{\infty}'/L_0'$, and hence that
$[L_{\infty}/L_0]=[L_{\infty}'/L_0']$.~\qed\medskip

     We define the depth of $g(X)\in\A(k)$ to be the degree
of the leading term of ${(g(X)-X)/X}$; the depth of $g(X)=X$ is
taken to be $\infty$.  Let $\Gamma$ be a subgroup of $\A(k)$
which is isomorphic to $\Z_p$, and let $\gamma$ be a
generator for $\Gamma$.  For $n\ge0$ we define the $n$th
lower ramification break $i_n$ of $\Gamma$ to be the depth of
$\gamma^{p^n}$; this definition is independent of the choice
of $\gamma$.  The upper ramification breaks of $\Gamma$ are
defined by the formulas $b_0=i_0$ and
$b_n-b_{n-1}=(i_n-i_{n-1})/p^n$ for $n\ge1$.  The $b_n$ are
integers by Sen's theorem \cite{sen}.  It follows from
\cite[Cor.\:3.3.4]{WZp} that if
$\Gamma=\Phi([L_{\infty}/L_0])$ then $(i_n)_{n\ge0}$ and
$(b_n)_{n\ge0}$ are the lower and upper ramification
sequences of the $\Z_p$-extension $L_{\infty}/L_0$.  Note
that since $L_{\infty}/L_0$ is an arithmetically profinite
extension \cite[\S1]{cn}, the Hasse-Herbrand functions
$\phi_{L_{\infty}/L_0}$ and $\psi_{L_{\infty}/L_0}$ are
defined, and the lower and upper ramification breaks of
$L_{\infty}/L_0$ are related by the formulas
$b_n=\phi_{L_{\infty}/L_0}(i_n)$ and
$i_n=\psi_{L_{\infty}/L_0}(b_n)$ for $n\ge0$.  For future
use we recall the following facts about the upper
ramification breaks of a cyclic extension
(see for instance \cite[p.\:280]{mar}).

\begin{lemma} \label{breaks}
Let $K$ be a local field with residue characteristic $p$ and
let $L/K$ be a totally ramified $(\Z/p^n\Z)$-extension.  Let
$1\le e\le\infty$ be the $K$-valuation of $p$ and let
${b_0<b_1<\dots<b_{n-1}}$ be the upper ramification
breaks of $L/K$.  Then for $0\le i\le n-2$ we have: \\[.2cm]
(a) $1\le b_0\le pe/(p-1)$; \\[.2cm]
(b) If $b_i\le e/(p-1)$ then $pb_i\le b_{i+1}\le pe/(p-1)$;
\\[.2cm]
(c) If $b_i\ge e/(p-1)$ then $b_{i+1}=b_i+e$.
\end{lemma}

     Let $\z_0(k)$ denote the subspace of $\z(k)$ consisting
of $k$-isomorphism classes of $\Z_p$-extensions in
characteristic $0$ and let $\G_0(k)=\Phi(\z_0(k))$.

\begin{cor} \label{Z0}
$\Phi|_{\z_0(k)}:\z_0(k)\ra\G_0(k)$ is a homeomorphism.
\end{cor}

\proof For $1\le e<\infty$ let $\z_0^e(k)$ denote the subspace of
$\z_0(k)$ consisting of $k$-isomorphism classes of
$\Z_p$-extensions
$[L_{\infty}/L_0]$ such that the absolute ramification index
of $L_0$ is $e$, and let $\G_0^e(k)=\Phi(\z_0^e(k))$.  Then
$\z_0^e(k)$ is open in $\z(k)$.  It follows from Krasner's Lemma
that there are only finitely many isomorphism classes of local
fields $L_0$ with residue field $k$ and absolute ramification
index $e$.  For each such $L_0$ consider the set $\H_{L_0}$
of continuous homomorphisms $\chi:L_0^{\times}\ra\Z_p$ such
that $\chi(\O_{L_0}^{\times})=\Z_p$.  This set is compact,
and class field theory gives a continuous map from $\H_{L_0}$
onto the set of all elements of $\z_0^e(k)$ of the form
$[L_{\infty}/L_0]$.  Therefore $\z_0^e(k)$ is compact.  Since
$\Phi$ is a continuous bijection, it follows that
\begin{equation} \label{homeo}
\Phi|_{\z_0^e(k)}:\z_0^e(k)\lra\G_0^e(k)
\end{equation}
is a homeomorphism.

     Let $[\Gamma]\in\G_0^e(k)$ and let
$[L_{\infty}/L_0]\in\z_0^e(k)$ be such that
$\Phi([L_{\infty}/L_0])=[\Gamma]$.  Let
$(b_n)_{n\ge0}$ be the upper ramification sequence of
$L_{\infty}/L_0$ and $\Gamma$.  It follows from
Lemma~\ref{breaks} that there
is $M\ge1$ such that $b_n-b_{n-1}=e$ for all $n\ge M$.  If
$[\Gamma']\in\G_0(k)$ is sufficiently close to
$[\Gamma]$ then the first $M+2$ upper ramification breaks
$b_0',b_1',\dots,b_{M+1}'$ of
$\Gamma'$ are the same as those of $\Gamma$.
In particular, we have $b_M'-b_{M-1}'=b_{M+1}'-b_M'=e$.  Let
$[L_{\infty}'/L_0']$ be the unique element of $\z_0(k)$ such that
$\Phi([L_{\infty}'/L_0'])=[\Gamma']$.  Then the upper
ramification breaks of $L_{\infty}'/L_0'$ are the same as
those of $\Gamma'$, so by Lemma~\ref{breaks} the
absolute ramification index of $L_0'$ is $e$.  It follows
that $[\Gamma']\in\G_0^e(k)$, and hence that $\G_0^e(k)$
is open in $\G_0(k)$.  Since (\ref{homeo}) is a homeomorphism
for $1\le e<\infty$, we conclude
that $\Phi$ induces a homeomorphism between $\z_0(k)$ and
$\G_0(k)$.~\qed \medskip

     Let $\z_p(k)$ denote the subspace of $\z(k)$
consisting of $k$-isomorphism classes of $\Z_p$-extensions
in characteristic $p$ and let $\G_p(k)=\Phi(\z_p(k))$.  Using
\cite[3.3]{WZp} and
Krasner's Lemma one can show that $\Phi$ induces a
homeomorphism between $\z_p(k)$ and $\G_p(k)$.  But $\Phi$
itself is not a homeomorphism.  Indeed, if
$[L_{\infty}/L_0]\in\z_p(k)$ then by the methods of the next
section one can construct $[L_{\infty}'/L_0']\in\z_0(k)$ such
that $\Phi([L_{\infty}'/L_0'])$ is arbitrarily close to
$\Phi([L_{\infty}/L_0])$.  Since the distance between
$[L_{\infty}/L_0]$ and $[L_{\infty}'/L_0']$ is always 2, this
implies that $\Phi$ is not a homeomorphism.

     Since $\G_0^e(k)$ is compact, it follows from
Corollary~\ref{Z0} that the map
\begin{equation}
\Phi^{-1}|_{\G_0^e(k)}:\G_0^e(k)\lra\z_0^e(k)
\end{equation}
is uniformly continuous.  From this fact we deduce the following
non-effective version of Theorem~\ref{main}:

\begin{cor}
Let $e\ge1$ and let $k$ be a finite field of characteristic
$p$.  Then there is a nondecreasing function
$s:{\mathbb N}\ra{\mathbb N}\cup\{0\}$ such that
$\lim_{n\ra\infty}s(n)=\infty$ which has the following
property: Let $L_0,L_0'$ be finite extensions of $\Q_p$ with
residue field $k$ and absolute ramification index $e$, and let
$L_{\infty}/L_0$ and $L_{\infty}'/L_0'$ be totally ramified
$\Z_p$-extensions such that
\begin{equation}
\Gamma_{L_{\infty}/L_0}^{(\pi_n)}\equiv
\Gamma_{L_{\infty}'/L_0'}^{(\pi_n')}\pmod{X^{m+1}}
\end{equation}
for some $m\ge1$ and some compatible sequences of uniformizers
$(\pi_n)$ for
$L_{\infty}/L_0$ and $(\pi_n')$ for $L_{\infty}'/L_0'$.  Then
$L_{s(m)}/L_0\cong_k L_{s(m)}'/L_0'$.
\end{cor}

\section{Truncated valuation rings} \label{tlr}

     In this section we give an overview of Deligne's theory
of extensions of truncated valuation rings.  For more details
see \cite{D}.

     Define a category $\T$ whose objects are triples
$(A,M,\epsilon)$ such that:
\begin{enumerate}
\item $A$ is an Artin local ring whose maximal ideal $m_A$
is principal and whose residue field is finite.
\item $M$ is a free $A$-module of rank 1.
\item $\epsilon:M\ra A$ is an $A$-module homomorphism whose
image is $m_A$.
\end{enumerate}
Let $S_1=(A_1,M_1,\epsilon_1)$
and $S_2=(A_2,M_2,\epsilon_2)$ be elements of $\T$.  A
morphism from $S_1$ to $S_2$ is a triple $f=(r,\mu,\eta)$,
where $r$ is a positive integer, $\mu:A_1\ra A_2$ is a ring
homomorphism, and $\eta:M_1\ra M_2^{\otimes r}$ is an
$A_1$-module homomorphism.  These must satisfy
$\mu\circ\epsilon_1=\epsilon_2^{\otimes r}\circ\eta$, and the
map $M_1\otimes_{A_1}A_2\ra M_2^{\otimes r}$ induced by $\eta$
must be an isomorphism of $A_2$-modules.  Let
$S_3=(A_3,M_3,\epsilon_3)$ be another element of $\T$, and
let $g=(s,\nu,\theta):S_2\ra S_3$ be a morphism.  Then the
composition of $g$ with $f$ is defined to be $g\circ f=
(sr,\nu\circ\mu,\theta^{\otimes r}\circ\eta)$.  Thus the
identity morphism on $S_1$ is
$(1,\mathrm{id}_{A_1},\mathrm{id}_{M_1})$,
and $f$ is an isomorphism if and only if $r=1$, $\mu$ is an
isomorphism, and $\eta$ is an isomorphism.
Let $f=(r,\mu,\eta)$ and
$f'=(r',\mu',\eta')$ be morphisms from $S_1$ to $S_2$,
and let $c$ be a positive integer.  We say
$f$ and $f'$ are $R(c)$-equivalent, or $f\equiv f'\pmod{R(c)}$,
if $r=r'$, $\mu$ and $\mu'$ induce the same map on residue
fields, and $\eta(x)-\eta'(x)\in m_{A_2}^{rc}M_2^{\otimes r}$
for all $x\in M_1$. 

     Let $f=(r,\mu,\eta):S_1\ra S_2$ be a $\T$-morphism.
We say that $(S_2,f)$ is an extension of $S_1$ if
$\length(A_2)=r\cdot\length(A_1)$.  We will often denote the
extension $(S_2,f)$ by $S_2/S_1$.
Let $(S_2,f)$ and $(S_3,g)$ be extensions of $S_1$.  A
morphism from $(S_2,f)$ to $(S_3,g)$ is defined to be a
$\T$-morphism $h:S_2\ra S_3$ such that $h\circ f=g$.
If $(S_2',f')$ is
an extension of $S_1'$, we say that $S_2'/S_1'$ is isomorphic
to $S_2/S_1$ if there are isomorphisms $i:S_1'\ra S_1$ and
$j:S_2'\ra S_2$ such that $j\circ f'=f\circ i$.

     Let $K$ be a local field and let $e$ be a
positive integer.  Define the $e$-truncation $\Tr_e(K)$ of
$K$ to be the triple $(A,M,\epsilon)$ consisting of the ring
$A=\O_K/\P_K^e$, the $A$-module $M=\P_K/\P_K^{e+1}$,
and the $A$-module homomorphism $\epsilon:M\ra A$ induced
by the inclusion $\P_K\hookrightarrow\O_K$.  It is clear that
$\Tr_e(K)$ is an element of $\T$.  Conversely, every element
of $\T$ is isomorphic to $\Tr_e(K)$ for some finite extension
$K$ of $\Q_p$ and some $e\ge1$ (cf.~\cite[1.2]{D}).

     Let $K$ and $L$ be local fields, let $\sigma:K\ra L$ be
an embedding, and let $r$ be the ramification
index of $L$ over $\sigma(K)$.  Define a morphism
\begin{equation}
f_{\sigma}=(r,\mu_{\sigma},\eta_{\sigma}):
\Tr_e(K)\lra\Tr_{re}(L)
\end{equation}
where
\begin{align}
\mu_{\sigma}&:\O_K/\P_K^e\lra\O_L/\P_L^{re} \\
\eta_{\sigma}&:\P_K/\P_K^{e+1}\lra\P_L^{r}/\P_L^{re+r}\cong
(\P_L/\P_L^{re+1})^{\otimes r}
\end{align}
are induced by $\sigma$.  Then $(\Tr_{re}(L),f_{\sigma})$ is
an extension of $\Tr_e(K)$.  If $L$ is a finite extension of
$K$ with ramification index $r$ we write
$f_{L/K}=(r,\mu_{L/K},\eta_{L/K})$ for the morphism
from $\Tr_e(K)$ to $\Tr_{re}(L)$ induced by the inclusion
$K\hookrightarrow L$.  The following proposition shows that
all extensions of $\Tr_e(K)$ are produced by this construction.

\begin{prop} \label{fields}
(\cite[Lemme 1.4.4]{D}) Let $K$ be a local field,
let $e\ge1$, and let $(T,f)$ be an extension of
$\Tr_e(K)$, with
$f=(r,\mu,\eta)$.  Then there is a finite extension
$L/K$ such that $(T,f)\cong(\Tr_{re}(L),f_{L/K})$.
\end{prop}

     Let $d\ge0$ be real, let $L/K$ be a finite extension of
local fields, and let $N/K$ be the normal closure of $L/K$
in $L^{sep}$.  We denote the largest upper ramification break
of $L/K$ by $u_{L/K}$.  We say that $L/K$ satisfies condition
$C^d$ if $d>u_{L/K}$, or equivalently, if the ramification
subgroup $\Gal(N/K)^d$ is
trivial.  Let $\ext(K)^d$ denote the category whose objects
are finite extensions of $K$ which satisfy condition
$C^d$, and whose morphisms are $K$-inclusions.

     Let $S\in\T$ and let $(T,f)$ be an extension of $S$.
Then there are positive integers $r,e$ and
a finite extension of local fields $L/K$ such that
$T/S\cong\Tr_{re}(L)/\Tr_e(K)$.  Given $0\le d\le e$
we say that $T/S$ satisfies condition $C^d$ if $L/K$
satisfies condition $C^d$.  This definition is independent
of the choice of $L/K$.  One can associate ramification
data to the extension
$T/S$.  In particular, the Hasse-Herbrand functions
$\phi_{T/S}$ and $\psi_{T/S}$ are defined.  It follows from
\cite[1.5.3]{D} that if $T/S$ satisfies condition $C^e$
then $\phi_{T/S}=\phi_{L/K}$ and $\psi_{T/S}=\psi_{L/K}$.

     Let $S\in\T$.  We define a category $\ext(S)^d$
whose objects are extensions of $S$ which satisfy condition
$C^d$.  An $\ext(S)^d$-morphism from $(T_1,f_1)$ to $(T_2,f_2)$
is defined to be an $R(\psi_{T_1/S}(d))$-equivalence class
of morphisms from $(T_1,f_1)$ to $(T_2,f_2)$.  The main result
of \cite{D} is the following.

\begin{theorem} \label{equiv}
(\cite[Th\'eor\`eme~2.8]{D})  Let $K$ be a local field and
let $e$ be a positive integer.
Then the functor from $\ext(K)^e$ to $\ext(\Tr_e(K))^e$ which
maps $L/K$ to $\Tr_{re}(L)/\Tr_e(K)$ is an equivalence of
categories.
\end{theorem}

     The proof of Theorem~\ref{main} depends on the following
application of Theorem~\ref{equiv}:

\begin{cor} \label{induced}
Let $e$ be a positive integer and let $L/K$ and $L'/K$ be finite
extensions of local fields which have ramification index
$r$ and satisfy condition $C^e$.  Let $\tau\in\Aut(K)$
and let $j:\Tr_{re}(L')\ra\Tr_{re}(L)$ be an isomorphism such
that $j\circ f_{L'/K}=f_{L/K}\circ f_{\tau}$.  Then there
is a unique isomorphism $\gamma:L'\ra L$ such that
$j\equiv f_{\gamma}\pmod{R(\psi_{L/K}(e))}$ and $\gamma|_K=\tau$.
\end{cor}

\proof Let $i:K\hookrightarrow L$ and $i':K\hookrightarrow L'$
be the inclusion maps.  Then $(\Tr_{re}(L'),f_{L'/K})$ and
$(\Tr_{re}(L),f_{L/K}\circ f_{\tau})$ are elements of
$\ext(\Tr_e(K))^e$ which are induced by $i'$
and $i\circ\tau$.  Since $j$ gives an
isomorphism between these extensions, by Theorem~\ref{equiv}
there is a unique isomorphism $\gamma:L'\ra L$ such that
$j\equiv f_{\gamma}\pmod{R(\psi_{L/K}(e))}$ and
$\gamma\circ i'=i\circ\tau$. \qed

\section{Recognizing cyclotomic extensions}
\label{appear}

     Before proving Theorem~\ref{main} we prove the
following result, which may be viewed as a special case of
the theorem.  An analogous result in the setting of the field
of norms is proved in \cite[Prop.\:3]{WZp}.

\begin{prop} \label{root}
Let $p>2$ and let $F/\Q_p$ be a finite tamely ramified
extension with ramification index $e$.  Set
$s=(p-1)/\gcd(e,p-1)$ and $e_0=e/\gcd(e,p-1)$.  Let $m\ge1$,
let $E/F$ be a totally ramified cyclic extension of degree
$sp^m$, and let $d$ be an integer such that $p\nmid d$ and
the image of $d$ in $(\Z/p^{m+1}\Z)^{\times}$ has order $sp^m$.
Assume there is $\alpha\in E$ such that $v_{E}(\alpha-1)=e_0$
and a generator $\tau$ for $\Gal(E/F)$ such that
$\tau(\alpha)\equiv\alpha^d\pmod{\P_E^n}$ for some
$n>e_0p^m$.  Then there is a primitive $p^{m+1}$th root of
unity $\xi\in\Q_p^{alg}$ such that $v_E(\alpha-\xi)\geq
(gs+e_0)p^m$, where
\begin{equation} \label{gdef}
g=\left\lceil\frac{n-e_0(p^{m+1}+p^m-1)}{sp^m}\right\rceil.
\end{equation}
In particular, if $n>e_0(p^{m+1}+p^m-1)$ then $E=F(\xi)$.
\end{prop}

     For $t\in\Z$ let $f(t)$ denote the maximum value of
$v_E(\tau(\beta)\beta^{-1}-d)$ as $\beta$ ranges over the
compact set $C_t=\{\beta\in E:v_E(\beta)=t\}$.  Since
$\N_{E/F}(d)\not=1$ we have $\tau(\beta)\beta^{-1}\not=d$
for all $\beta\in C_t$, so $f(t)\in\Z$.  The proof of
Proposition~\ref{root} depends on the following lemma:

\begin{lemma} \label{shift}
Let $t\in\Z$ and set $t_0=t-e_0p^m$.  Then
\begin{equation} \label{ft}
f(t)=
\begin{cases}
0&\text{if $s\nmid t_0$,} \\
e_0(p^{v_p(t_0)+1}-1)&\text{if $s\mid t_0$ and $v_p(t_0)<m$,} \\
e_0(p^{m+1}-1)&\text{if $s\mid t_0$ and $v_p(t_0)\ge m$.}
\end{cases}
\end{equation}
\end{lemma}

     The proof of Lemma~\ref{shift} uses the following result,
which follows easily from Proposition~8 in \cite[V]{cl}.

\begin{prop} \label{norms}
Let $E$ be a local field and let
$M/E$ be a finite totally ramified Galois extension.
Let $d$ be a positive
integer and let $x,y$ be elements of $\O_M^{\times}$ such that
$x\equiv y\pmod{\P_M^{\psi_{M/E}(d)+1}}$.  Then
$\N_{M/E}(x)\equiv\N_{M/E}(y)\pmod{\P_E^{d+1}}$.
\end{prop}

\noindent {\em Proof of Lemma~\ref{shift}:}
Since $E/F$ is totally ramified and $p\nmid e_0$
we can write $\alpha=1+c\pi_E^{e_0}$, where $c\in\O_F^{\times}$
and $\pi_E$ is a uniformizer for $E$.  Since $v_p(d^s-1)=1$
and $v_E(\alpha^p-1)=pe_0$, we have
$v_E(\tau^s(\alpha)-\alpha)=pe_0$.  It follows that
$v_E(\tau^s(\pi_E^{e_0})-\pi_E^{e_0})=pe_0$, and hence that
$v_E(\tau^s(\pi_E^{e_0})\pi_E^{-e_0}-1)=(p-1)e_0$.
Since $\tau^s$ has order $p^m$ we have
$v_E(\tau^s(\pi_E)\pi_E^{-1}-1)\ge1$, and hence
$v_E(\tau^s(\pi_E)\pi_E^{-1}-1)=(p-1)e_0$.  Let $T/F$ denote
the maximum tamely ramified subextension of $E/F$.  Then $T$
is the subfield of $E$ fixed by $\langle\tau^s\rangle$, so
the smallest (upper and lower) ramification break of $E/T$ is
$(p-1)e_0$.  Since $(p-1)e_0=v_T(p)$, by
Lemma~\ref{breaks}(c) we deduce that for $0\le i<m$ the
$i$th upper ramification break of $E/T$ is $(p-1)e_0(i+1)$.
It follows that the $i$th lower ramification break of $E/T$ is
$\psi_{E/T}((p-1)e_0(i+1))=e_0(p^{i+1}-1)$.

     Let $\gamma=\alpha-1=c\pi_E^{e_0}$.  Then
\begin{equation} \label{quot}
\frac{\tau(\gamma)}{\gamma}\equiv\frac{(1+\gamma)^d-1}{\gamma}
\equiv d\pmod{\P_E^{e_0}}.
\end{equation}
Since the smallest upper ramification break of $E/T$ is
$(p-1)e_0$, we have $\psi_{E/T}(e_0-1)=e_0-1$.  Applying
Proposition~\ref{norms} to (\ref{quot}) we get
\begin{equation} \label{Nu}
\N_{E/T}\left(\frac{\tau(\gamma)}{\gamma}\right)\equiv
\N_{E/T}(d)\pmod{\P_T^{e_0}}.
\end{equation}
Let $\delta=\N_{E/T}(\gamma)$.
Since $\Gal(E/F)$ is commutative, (\ref{Nu}) reduces to
\begin{equation} \label{rhopi}
\frac{\tau(\delta)}{\delta}\equiv d^{p^m}\equiv d
\pmod{\P_T^{e_0}}.
\end{equation}

     Since both sides of (\ref{ft})
depend only on the congruence class of $t$ modulo $sp^m$,
we may assume $e_0p^m\le t<(e_0+s)p^m$.
Let $\beta$ be an element of $E$ such that $v_E(\beta)=t$,
and set $\kappa=\beta\delta^{-1}$.  Then by (\ref{rhopi}) we
have
\begin{align}
\frac{\tau(\beta)}{\beta}-d&=\frac{\tau(\kappa)}{\kappa}\cdot
\frac{\tau(\delta)}{\delta}-d \\
&\equiv \left(\frac{\tau(\kappa)}{\kappa}-1\right)d
\pmod{\P_E^{e_0p^m}}. \label{rhonu}
\end{align}
For $0<t_0<sp^m$ let $g(t_0)$ denote
the maximum value of $v_E(\tau(\kappa)\kappa^{-1}-1)$ as
$\kappa$ ranges over $C_{t_0}$.
It follows from Sen's argument in \cite[p.\:35]{sen}
that $g(t_0)$ is equal to the ramification number
$v_E(\tau^{t_0}(\pi_E)\pi_E^{-1}-1)$ of $\tau^{t_0}$.
Thus if $s\nmid t_0$ then $g(t_0)=0$, while if $s\mid t_0$
and $v_p(t_0)=i$ then $g(t_0)=e_0(p^{i+1}-1)$ is the $i$th
lower ramification break of $E/T$.  It follows that
$g(t_0)<e_0p^m$ for $0<t_0<sp^m$.  Hence by (\ref{rhonu})
we get $g(t_0)=f(t_0+e_0p^m)$.  This proves
the lemma for all $t$ such that $e_0p^m<t<(e_0+s)p^m$.

     It remains to prove the lemma for $t=e_0p^m$.
For $e_0p^m\le t<(e_0+s)p^m$ let $\beta_t$ be an element
of $E$ such that $v_E(\beta_t)=t$ and $v_E((\tau-d)\beta_t)$
is maximized.  Let $\Lambda$ denote the $\O_F$-lattice spanned
by the $\beta_t$.  Then $\Lambda=\pi_E^{e_0p^m}\O_E$ is an
ideal in $\O_E$, and hence
$(\tau-d)\Lambda$ is contained in $\Lambda$.
It follows from the maximality of $v_E((\tau-d)\beta_t)$ that
the integers $v_E((\tau-d)\beta_t)$ for
$e_0p^m\le t<(e_0+s)p^m$
represent distinct congruence classes modulo $sp^m$.
Therefore $\Lambda/(\tau-d)\Lambda$ is an $\O_F$-module of
length
\begin{equation}
\sum_{t=e_0p^m}^{(e_0+s)p^m-1}v_E((\tau-d)\beta_t)
\;\;-\sum_{t=e_0p^m}^{(e_0+s)p^m-1}v_E(\beta_t)=
\sum_{t=e_0p^m}^{(e_0+s)p^m-1}f(t).
\end{equation}
It follows that
\begin{equation}
\sum_{t=e_0p^m}^{(e_0+s)p^m-1}f(t)=v_E(\det(\tau-d)).
\end{equation}
Since the characteristic polynomial of the $F$-linear map
$\tau:E\ra E$ is $h(X)=X^{sp^m}-1$, the determinant of
$\tau-d$ is $\pm h(d)=\pm(d^{sp^m}-1)$.  Therefore we have
\begin{align}
\sum_{t=e_0p^m}^{(e_0+s)p^m-1}f(t)&=v_E(d^{sp^m}-1) \\
&=(m+1)e_0(p^{m+1}-p^m).
\end{align}
Solving for $f(e_0p^m)$ in terms of the known values of
$f(t)$ gives $f(e_0p^m)=e_0(p^{m+1}-1)$, which completes
the proof of the lemma. \qed \medskip

\noindent
{\em Proof of Proposition~\ref{root}:} It follows from the
hypotheses that $v_E(\alpha^{p^{m+1}}-1)\geq e_0p^{m+1}$, and that
\begin{equation}
v_E(\tau(\alpha^{p^{m+1}})-\alpha^{dp^{m+1}})\ge
n+(m+1)e_0(p^{m+1}-p^m).
\end{equation}
Let $\lambda=\log(\alpha^{p^{m+1}})$.  Then we have
\begin{equation}
v_E(\tau(\lambda)-d\lambda)=
v_E(\tau(\alpha^{p^{m+1}})-\alpha^{dp^{m+1}}),
\end{equation}
and hence
\begin{equation} \label{lsum}
v_E(\lambda)+v_E\!\left(\frac{\tau(\lambda)}{\lambda}-d\right)\ge
n+(m+1)e_0(p^{m+1}-p^m).
\end{equation}
Set $t=v_E(\lambda)=v_E(\alpha^{p^{m+1}}-1)$.  Then by
(\ref{lsum}) we get
\begin{equation} \label{tft}
t+f(t)\ge n+(m+1)e_0(p^{m+1}-p^m),
\end{equation}
where $f(t)$ is the function defined in Lemma~\ref{shift}.

     If $f(t)>0$ then by Lemma~\ref{shift} we have
$t=e_0p^m+csp^i$ and $f(t)=e_0(p^{i+1}-1)$ 
for some $0\le i\le m$ and $c\in\Z$.
It follows from (\ref{tft}) that 
\begin{equation}
e_0p^m+csp^i+e_0(p^{i+1}-1)\ge n+(m+1)e_0(p^{m+1}-p^m),
\end{equation}
which implies
\begin{equation}
csp^i\ge e_0(p^{m+1}-p^{i+1})+(m+1)e_0(p^{m+1}-p^m)
+n-e_0(p^{m+1}+p^m-1).
\end{equation}
Dividing by $sp^i$ and using the fact that $s$ divides
$p-1$ we get
\begin{equation}
c\ge e_0\frac{p^{m+1}-p^{i+1}}{sp^i}+(m+1)e_0\frac{p^{m+1}-p^m}{sp^i}
+\left\lceil\frac{n-e_0(p^{m+1}+p^m-1)}{sp^i}\right\rceil.
\end{equation}
It follows that
\begin{equation} \label{rbound}
t\ge e_0p^m+e_0(p^{m+1}-p^{i+1})+(m+1)e_0(p^{m+1}-p^m)
+sp^i\left\lceil\frac{n-e_0(p^{m+1}+p^m-1)}{sp^i}\right\rceil.
\end{equation}
The minimum value of right hand side of (\ref{rbound}) for
$0\le i\le m$ is achieved when $i=m$.  Therefore the
inequality
\begin{equation} \label{tbound}
t\ge e_0p^m+(m+1)e_0(p^{m+1}-p^m)+
sp^m\left\lceil\frac{n-e_0(p^{m+1}+p^m-1)}{sp^m}\right\rceil
\end{equation}
holds for all $t$ such that $f(t)>0$.  If $f(t)=0$ then by
(\ref{tft}) we have
\begin{equation}
t\ge n+(m+1)e_0(p^{m+1}-p^m),
\end{equation}
which implies that (\ref{tbound}) holds
in this case as well.  Thus (\ref{tbound}) is valid in
general.

     Let $\zeta\in\Q_p^{alg}$ be a primitive $p^{m+1}$th root
of unity, and choose $0\leq j<p^{m+1}$ to maximize
$w=v_E(\alpha-\zeta^j)$.  For $0\leq i<p^{m+1}$ we have
\begin{equation}
v_E(\alpha-\zeta^i) \geq \min\{w,v_E(\zeta^j-\zeta^i)\},
\end{equation}
with equality if $w>v_E(\zeta^j-\zeta^i)$.  Since $w\geq
v_E(\alpha-\zeta^i)$, this implies that for $i\not=j$ we have
$v_E(\alpha-\zeta^i)\leq v_E(\zeta^j-\zeta^i)=e_0p^{v_p(i-j)}$.
Since
\begin{equation}
\alpha^{p^{m+1}}-1 = (\alpha-1)(\alpha-\zeta)(\alpha-\zeta^2)\dots
(\alpha-\zeta^{p^{m+1}-1}),
\end{equation}
by defining $p^{v_p(0)}=0$ we get
\begin{align}
t=v_E(\alpha^{p^{m+1}}-1)&\leq w+\sum_{i=0}^{p^{m+1}-1}\,e_0p^{v_p(i-j)} \\
&=w+(m+1)e_0(p^{m+1}-p^m). \label{m}
\end{align}
By comparing (\ref{tbound}) with (\ref{m}) we conclude that
$w\ge gsp^m+e_0p^m$, where $g$ is the integer defined in
(\ref{gdef}).  Set $\xi=\zeta^j$; then
$v_E(\alpha-\xi)=w\ge gsp^m+e_0p^m$.  Since
$v_E(\alpha-\xi)>e_0$, we have $v_E(\xi-1)=v_E(\alpha-1)=e_0$,
so $\xi$ is a primitive $p^{m+1}$th root of unity.
If $n>e_0(p^{m+1}+p^m-1)$ then $g\ge1$, and hence
$v_E(\alpha-\xi)> e_0p^m$.  Therefore by
Krasner's Lemma we have $F(\alpha)\supset F(\xi)$.  Since
$E\supset F(\alpha)$ and $[F(\xi):F]\ge sp^m=[E:F]$, this
implies $E=F(\xi)$. \qed 

\section{Proof of Theorem~\ref{main}} \label{proof}

     In this section we prove Theorem~\ref{main} in a
somewhat generalized form.  Let $1\le a\le ep^n$ and
$1\le m\le n$.  We will show that $[\omega(L)\cap L':K]\ge p^m$
whenever
\begin{equation} \label{hyp}
h_{\pi_L}^{\sigma}(X)\equiv h_{\pi_{L'}}^{\sigma'}(X)\pmod{X^a}
\end{equation}
and $a$ and $m$ satisfy certain inequalities, which are
specified in Theorem~\ref{conditions}.
We then show in Lemmas \ref{cMEq} and
\ref{psid} that the values $a=ep^n$ and $m=m_0$ given
in Theorem~\ref{main} satisfy these inequalities.  To
motivate the proof we first prove an analog of
Theorem~\ref{main} for local fields of characteristic $p$.

\begin{prop}
Let $K$ be a local field of characteristic $p$ with
residue field $k$ and let $L/K$, $L'/K$ be totally ramified
$(\Z/p^n\Z)$-extensions.  Let $u_{L/K}$ be the largest
upper ramification break of $L/K$, let $e>u_{L/K}$, and let
$h(X)\in k[X]$.  Assume there exist uniformizers
$\pi_L,\pi_{L'}$ for $L,L'$ and generators $\sigma,\sigma'$
for $\Gal(L/K)$, $\Gal(L'/K)$ such that
\begin{alignat}{2}
\sigma(\pi_L)&\equiv\pi_Lh(\pi_L)&&\pmod{\P_L^{ep^n+1}}
\label{spi} \\
\sigma'(\pi_{L'})&\equiv\pi_{L'}h(\pi_{L'})&&
\pmod{\P_{L'}^{ep^n+1}}. \label{spibis}
\end{alignat}
Then the extensions $L/K$ and $L'/K$ are $k$-isomorphic.
\end{prop}

\proof Let $\alpha:L'\ra L$ be the unique $k$-isomorphism
such that $\alpha(\pi_{L'})=\pi_L$.
Since $\pi_K=\N_{L/K}(\pi_L)$ and
$\pi_K'=\N_{L'/K}(\pi_{L'})$ are
uniformizers for $K$ there is a unique $k$-automorphism
$\tau$ of $K$ such that $\tau(\pi_K')=\pi_K$.
It follows from (\ref{spi}) and (\ref{spibis}) that
\begin{equation}
\tau(\pi_K')\equiv\alpha(\pi_K')\pmod{\P_L^{ep^n+1}}.
\end{equation}
Therefore the induced maps
\begin{align}
f_{L/K}&:\Tr_e(K)\lra\Tr_{ep^n}(L) \\
f_{L'/K}&:\Tr_e(K)\lra\Tr_{ep^n}(L') \\
f_{\tau}&:\Tr_e(K)\lra\Tr_e(K) \\
f_{\alpha}&:\Tr_{ep^n}(L')\lra\Tr_{ep^n}(L)
\end{align}
satisfy $f_{\alpha}\circ f_{L'/K}=f_{L/K}\circ f_{\tau}$.
Since $u_{L/K}=u_{L'/K}<e$, both $L/K$ and $L'/K$ satisfy
condition $C^e$.  Therefore by Corollary~\ref{induced} there is
a $k$-isomorphism $\gamma:L'\ra L$ such that $\gamma|_K=\tau$.
Hence $L'/K\cong_kL/K$. \qed \medskip

     To apply this method in characteristic 0 we replace
the fields $K,L,L'$ with cyclotomic extensions.  This makes
our fields resemble local fields of characteristic $p$ and
allows us to replace (\ref{hyp}) with a congruence modulo a
higher power of $X$.  Let $\zeta\in\Q_p^{alg}$ be a primitive
$p^{m+1}$th root of unity and set $M=L(\zeta)$.  Then $M/K$
is an abelian extension whose Galois group may be identified
with a subgroup of ${\Gal(L/K)\times\Gal(K(\zeta)/K)}$.
We will use the theory of
truncated local rings outlined in Section~\ref{tlr} to define an
extension $M'/L'$ which corresponds to $M/L$.  We will
then use Proposition~\ref{root} to show that in fact
$M'=L'(\zeta)$.  Let $L_0/K$, $L_0'/K$ be the subextensions
of $L/K$, $L'/K$ of degree $p^m$.  Using Corollary~\ref{induced}
we will show that $L_0(\zeta)/K\cong_k L_0'(\zeta)/K$, from which
it will follow that $L_0/K\cong_k L_0'/K$.

\begin{figure}
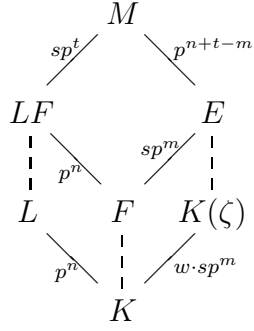

\[\begin{array}{rcccl}
&&M \\[.1cm]
&\ufs{\hspace{-.1cm}sp^t}&&\ubs{p^{n+t-m}} \\[.1cm]
LF\!\!\!&&&&E \\[.1cm]
\!\!\raisebox{4pt}{$\scriptstyle$}\dashline\,&\:\lbs{p^n}&&\ufs{\hspace{-.2cm}sp^m}&
\;\dashline\,\raisebox{4pt}{$\scriptstyle$}\! \\[.1cm]
L&&F&&\hspace{-.3cm}K(\zeta) \\[.1cm]
&\lbs{p^n}&\hspace{-.3cm}\dashline\,\raisebox{4pt}{$\scriptstyle$}\!\!\!\!\!&\lfs{w\cdot sp^m} \\[.1cm]
&&K
\end{array}\]
\caption{Dashed lines represent unramified extensions.}
\label{diagram}
\end{figure}

     Let $w$ denote the residue class degree and $sp^m$ the
ramification index of $K(\zeta)/K$.  Then the
ramification index of $M/L$ is equal to $sp^t$
for some $0\le t\le m$.  Let $F/K$ be the maximum unramified
subextension of $M/K$ and let $E/K(\zeta)$
be the maximum unramified subextension of $M/K(\zeta)$.
Then $E/F$ is a totally ramified cyclic extension of degree
$sp^m$, and $M/E$ is a totally ramified cyclic extension of
degree $p^{n+t-m}$ (see Figure~\ref{diagram}).

     In order to state our generalized version of
Theorem~\ref{main} we must first compute the
ramification data of the extension $L/K$.  Let $y$ be the
smallest upper ramification break of $L/K$ which exceeds
$\frac{1}{p-1}\cdot e$; if all the upper ramification breaks
of $L/K$ are $\le\frac{1}{p-1}\cdot e$, let $y$ be the
largest upper ramification break of $L/K$.  By
Lemma~\ref{breaks}(b) we have $y\le(1+\frac{1}{p-1})e$.
Suppose that $y=b_h$, where $b_0<b_1<\dots<b_{n-1}$ are the
upper ramification breaks of $L/K$, and let $z=\psi_{L/K}(y)$
be the corresponding lower break.  It follows from
Lemma~\ref{breaks}(c) that for $h\le i<n$ the $i$th upper
ramification break of $L/K$ is $b_i=y+(i-h)e$.  Therefore
for $h\le i<n$ the $i$th lower ramification break of $L/K$ is
\begin{align}
\psi_{L/K}(y+(i-h)e)&=z+ep^{h+1}+\dots+ep^i \\
&=z+ep^{h+1}\cdot\frac{p^{i-h}-1}{p-1}. \label{lowerLK}
\end{align}
The largest upper ramification break $u_{L/K}=b_{n-1}$ of $L/K$
is equal to $y+(n-h-1)e$.  Since $y\le(1+\frac{1}{p-1})e$ this
implies $u_{L/K}\le(n-h+\frac{1}{p-1})e$.  It follows that
\begin{equation} \label{mbound}
\psi_{L/K}\left(\left(n-h+1+\frac{1}{p-1}\right)e\right)>ep^n.
\end{equation}
Thus if $m$ satisfies $\psi_{L/K}((m+1+\frac{1}{p-1})e)<ep^n$
then $m\le n-h-1$.

     Set $e_0=es/(p-1)$ and define
\begin{equation}
q=\begin{cases}
((y-e)s+e_0)p^m&\text{if $h=0$ and $y>e$}, \\
e_0p^m&\text{otherwise}.
\end{cases}
\end{equation}
Also set $r=q+e_0(p^{m+1}-1)$.  Note that the integers $t$,
$q$, $r$ and the fields $M$, $E$ all depend on $m$.

\begin{theorem} \label{conditions}
Let $p>3$ and let $K$ be a finite tamely ramified extension
of $\Q_p$ with ramification index $e$.  Let $L/K$ and $L'/K$
be totally ramified $(\Z/p^n\Z)$-extensions such that $L/K$
is contained in a $\Z_p$-extension $L_{\infty}/K$.  Let $1\le
a\le ep^n$ and assume that there are
generators $\sigma$, $\sigma'$ for $\Gal(L/K)$, $\Gal(L'/K)$
and uniformizers $\pi_L$, $\pi_{L'}$ for $L$, $L'$
such that
$h_{\pi_L}^{\sigma}(X)\equiv
h_{\pi_{L'}}^{\sigma'}(X)\pmod{X^a}$.
Suppose there exists $1\le m\le n$ such that
the following three inequalities are satisfied:
\begin{align}
\psi_{M/L}(a)&>[M:E]q=p^{n+t-m}q \label{cond1} \\[.2cm]
\psi_{M/L}(a)&>\psi_{M/E}(r)     \label{cond2} \\[.2cm]
\psi_{M/L}(a)&>\psi_{M/K}(u_{L/K}).    \label{cond3}
\end{align}
Then there is $\omega\in\Gal(\Q_p^{alg}/\Q_p)$ such
that $\omega(K)=K$, $\omega$ induces the identity on $k$,
and $[L\cap\omega(L'):K]\geq p^m$.
\end{theorem}

     The following lemmas will be used to compute and bound
the ramification breaks of the various extensions used in the
proof of Theorem~\ref{conditions}.

\begin{lemma} \label{tame}
Let $K$ be a finite extension of $\Q_p$, let $L/K$ be a
finite tamely ramified extension with ramification index $e$,
and let $M/L$ be a finite Galois extension which is not
tamely ramified.  Then the positive lower ramification
breaks of $M/L$ are the same as the positive lower
ramification breaks of $M/K$, and the positive upper
ramification breaks of $M/L$ are $e$ times the positive upper
ramification breaks of $M/K$.
\end{lemma}

\proof The positive lower ramification breaks of
$M/K$ are the values $x>0$ such that
$\phi_{M/K}(x)=\phi_{L/K}\circ\phi_{M/L}(x)$ is not
differentiable.  It follows from \cite[Prop.\:A.4.2]{D} that
for $x>0$ we have $\phi_{L/K}(x)=e^{-1}\cdot x$.  Therefore
the positive
lower ramification breaks of $M/K$ and $M/L$ are the same.
Let $l_0<l_1<\dots<l_{n-1}$ be the positive lower
ramification breaks of $M/K$ and $M/L$.  Then the $i$th
positive upper ramification break of $M/L$ is $\phi_{M/L}(l_i)$,
and the $i$th positive upper ramification break of $M/K$ is
$\phi_{M/K}(l_i)=\phi_{L/K}\circ\phi_{M/L}(l_i)=
e^{-1}\cdot\phi_{M/L}(l_i)$. \qed

\begin{lemma} \label{change}
Let $K$ be a finite extension of $\Q_p$ and let $L/K$ be a
finite cyclic extension whose positive upper ramification
breaks are $b_0<b_1<\dots<b_{n-1}$.  Let $E/K$ be a finite
Galois extension and write the ramification index of $LE/E$
in the form $up^r$ with $p\nmid u$.  Then the positive upper
ramification breaks $\beta_{n-r}<\beta_{n-r+1}<\dots<\beta_{n-1}$
of $LE/E$
satisfy $\beta_i\le\psi_{E/K}(b_i)$ for $n-r\le i<n$.
In particular, $u_{LE/E}\le\psi_{E/K}(u_{L/K})$.
\end{lemma}

\proof We first prove the second statement.  Let
$m=\psi_{E/K}(u_{L/K})$, and suppose that
$m<u_{LE/E}$.  Since $L/K$ and
$LE/E$ are abelian, $u_{L/K}$ and $u_{LE/E}$ are integers.
Therefore $m$ is also an integer.  It
follows that $m+1\le u_{LE/E}$, and hence that
$\Gal(LE/E)^{m+1}$ is nontrivial.  Since the reciprocity map
$\omega_{LE/E}:E^{\times}\ra\Gal(LE/E)$ maps $1+\P_E^{m+1}$
onto $\Gal(LE/E)^{m+1}$ (see for instance Corollary~3 to
Theorem~2 in
\cite[XV\,\S2]{cl}), there is $\alpha\in1+\P_E^{m+1}$ such
that $\sigma=\omega_{LE/E}(\alpha)$ is not the identity.
It follows from the functorial properties of the
reciprocity map \cite[XI\,\S3]{cl} that
$\sigma|_L=\omega_{L/K}(\N_{E/K}(\alpha))$.
Using Proposition~\ref{norms} we see that
$\N_{E/K}(\alpha)\in1+\P_K^{u_{L/K}+1}$.  Since
$\omega_{L/K}(1+\P_K^{u_{L/K}+1})=\Gal(L/K)^{u_{L/K}+1}$ is
trivial this implies $\sigma|_L=\id_L$.
Since the restriction map $\Gal(LE/E)\ra\Gal(L/K)$ is
one-to-one, this is a contradiction.  Therefore
$u_{LE/E}\le\psi_{E/K}(u_{L/K})$.

     To prove the first statement, for $0\le j\le r-1$ let
$L^j/K$ be the unique
subextension of $L/K$ such that $[L:L^j]=p^j$.  The
restriction map $\Gal(LE/E)\ra\Gal(L/K)$
induces an isomorphism between $\Gal(LE/L)$ and
$\Gal(L/(L\cap E))$.  Since $p^j$ divides the ramification
index of $LE/E$, we see that $L^j\supset L\cap E$, that
$L/L^j$ is totally ramified, and that $\Gal(LE/L^jE)$ is the
unique subgroup of $\Gal(LE/E)$ of order $p^j$.  It follows
that $u_{L^j/K}=b_{n-j-1}$ and $u_{L^jE/E}=\beta_{n-j-1}$.
Applying the second statement we get
$\beta_{n-j-1}\le\psi_{E/K}(b_{n-j-1})$ for $0\le j\le r-1$.
\qed \medskip

\noindent {\em Proof of Theorem \ref{conditions}:}
Since $K$, $L$, and $L'$ all have the same residue field
$k$, there is a unique $k$-linear ring isomorphism
$\mu:\O_{L'}/\P_{L'}^a\ra\O_L/\P_L^a$ such that
$\mu(\pi_{L'})\equiv\pi_L\pmod{\P_L^a}$, and a unique
$(\O_{L'}/\P_{L'}^a)$-module isomorphism
$\eta:\P_{L'}/\P_{L'}^{a+1}\ra\P_L/\P_L^{a+1}$ such that
$\eta(\pi_{L'})\equiv\pi_L\pmod{\P_L^{a+1}}$.  By combining
these isomorphisms we get an isomorphism
$i=(1,\mu,\eta)$ from $\Tr_a(L')$ to $\Tr_a(L)$.  Since
$h_{\pi_L}^{\sigma}(X)\equiv
h_{\pi_{L'}}^{\sigma'}(X)\pmod{X^a}$,
we have $i\circ f_{\sigma'}=f_{\sigma}\circ i$.

     Let $b=[M:LF]\cdot a$.  Then
$(\Tr_b(M),f_{M/L}\circ i)$ is an extension of $\Tr_a(L')$.
It follows from Proposition~\ref{fields} that this extension
comes from an extension of $L'$.  More precisely, there
is a finite extension $M'/L'$ and an isomorphism
\begin{equation} \label{j}
j=(1,\nu,\theta):\Tr_b(M')\lra\Tr_b(M)
\end{equation}
such that $j\circ f_{M'/L'}=f_{M/L}\circ i$.
Since $K(\zeta)/\Q_p(\zeta)$ is tamely ramified we have
$u_{K(\zeta)/\Q_p}=u_{\Q_p(\zeta)/\Q_p}=m$.  Using
Lemma~\ref{tame} we see that $u_{K(\zeta)/K}=e\cdot
u_{K(\zeta)/\Q_p}=me$.  It follows from Lemma~\ref{change} that
$u_{M/L}\le\psi_{L/K}(u_{K(\zeta)/K})=\psi_{L/K}(me)$.  By
assumption (\ref{cond2}) we have
\begin{equation}
\psi_{M/L}(a)>\psi_{M/E}(r)>\psi_{M/E}(e_0(p^m-1)).
\end{equation}
Since $\psi_{E/K}(me)=e_0(p^m-1)$ this implies
$\psi_{M/L}(a)>\psi_{M/K}(me)$.  Applying $\phi_{M/L}$
to this inequality gives $a>\psi_{L/K}(me)\ge u_{M/L}$.
Thus $M/L$, $\Tr_b(M)/\Tr_a(L)$,
$\Tr_b(M')/\Tr_a(L')$, and $M'/L'$ all satisfy
condition $C^a$.  It follows from Theorem~\ref{equiv} that
the field $M'$ is uniquely determined up to $L'$-isomorphism.
Let $c=\psi_{M/L}(a)$.  By Theorem~\ref{equiv} the
isomorphism $j$ in (\ref{j}) is uniquely determined up
to $R(c)$-equivalence.

\begin{lemma} \label{lift}
Let $\gamma\in\Gal(M/K)$ and let $t\in\Z$ be such that
$\gamma|_L=\sigma^t$.  Then there is a unique automorphism
$\gamma'$ of $M'$ such that $\gamma'|_{L'}={\sigma'}^t$ and
$j\circ f_{\gamma'}\equiv f_{\gamma}\circ j\pmod{R(c)}$.
The map $\gamma\mapsto\gamma'$ gives a faithful
$K$-linear action of $\Gal(M/K)$ on $M'$.
\end{lemma}

\proof For $\gamma\in\Gal(M/K)$ let
$f_{\gamma}'=j^{-1}\circ f_{\gamma}\circ j$ denote
the automorphism of $\Tr_b(M')$ induced by $f_{\gamma}$.
Using the identities
$j\circ f_{M'/L'}=f_{M/L}\circ i$,
$f_{\gamma}\circ f_{M/L}=f_{M/L}\circ f_{\sigma^t}$, and
$f_{\sigma^t}\circ i=i\circ f_{\sigma'^t}$ we find that
$f_{\gamma}'\circ f_{M'/L'}=f_{M'/L'}\circ f_{{\sigma'}^t}$.
Since $M'/L'$ satisfies condition $C^a$, by
Corollary~\ref{induced} there is a unique
$\gamma'\in\Aut(M')$ such that
$f_{\gamma'}\equiv f_{\gamma}' \pmod{R(c)}$ and
$\gamma'|_{L'}={\sigma'}^t$.  It follows that $j\circ
f_{\gamma'}\equiv f_{\gamma}\circ j\pmod{R(c)}$.
Since $\gamma'$ is uniquely
determined by $\gamma$ the map $\gamma\mapsto\gamma'$ is a
group homomorphism.  If $\gamma$ lies in the kernel of this
homomorphism then ${\sigma'}^t=1$, and hence $\sigma^t=1$.
Therefore $\gamma\in\Gal(M/L)$ and $\gamma$ induces the
identity on $\Tr_c(M)$.  Since $M/L$ satisfies
condition $C^a$ this implies $\gamma=1$. \qed \medskip

     It follows from this lemma that $M'/K$ is Galois, and
that the map
\begin{equation} \label{jhatdef}
\jhat:\Gal(M/K)\lra\Gal(M'/K)
\end{equation}
defined by $\jhat(\gamma)=\gamma'$ is an isomorphism.
Furthermore, for all $\gamma\in\Gal(M/K)$ we have
\begin{equation} \label{jhat}
f_{\gamma}\circ j\equiv j\circ f_{\jhat(\gamma)}\pmod{R(c)}.
\end{equation}
Since $M'$ is a Galois extension
of $L'$ which is uniquely determined up to $L'$-isomorphism,
$M'$ is uniquely determined as a subfield of $\Q_p^{alg}$.

\begin{lemma} \label{Zp}
Let $K$ be a finite extension of $\Q_p$ and let $L/K$ be a
$(\Z/p^n\Z)$-extension.  Then
$L$ is contained in a $\Z_p$-extension $L_{\infty}$ of $K$ if
and only if the group $\mu$ of $p$-power roots of unity in $K$
is contained in $\N_{L/K}(L^{\times})$.
\end{lemma}

\proof If $L$ is contained in a $\Z_p$-extension
$L_{\infty}$ of $K$ then there is a continuous homomorphism
$\chi:K^{\times}\ra\Gal(L_{\infty}/K)$ such that
$\chi(K^{\times})$ is dense in $\Gal(L_{\infty}/K)$ and
$\ker(\chi)\le\N_{L/K}(L^{\times})$.
It follows that $K^{\times}/\ker(\chi)$ has trivial torsion,
and hence that $\mu\le\ker(\chi)\le\N_{L/K}(L^{\times})$.  If
$\mu\le\N_{L/K}(L^{\times})$ then since
$\N_{L/K}(L^{\times})$ has index $p^n$ in $K^{\times}$, the
group $\tilde{\mu}$ of all roots of unity in $K$ is contained in
$\N_{L/K}(L^{\times})$.  Since $K^{\times}/\tilde{\mu}\cong
\Z\times\Z_p^{[K:\Q]}$ there is a closed subgroup
$H$ of $\N_{L/K}(L^{\times})$ such that $K^{\times}/H\cong\Z_p$.
Then $H$ corresponds by class field theory to a
$\Z_p$-extension $L_{\infty}$ of $K$ which contains $L$.~\qed
\medskip

     Since $L$ is contained in a $\Z_p$-extension
$L_{\infty}$ of $K$, the field $M=LE$ is contained in the
$\Z_p$-extension $L_{\infty}E$ of $E$.  Therefore by
Lemma~\ref{Zp} there is $\alpha\in\O_M^{\times}$ such that
$\N_{M/E}(\alpha)=\zeta$.  Let $\tau$ be an element of
$\Gal(M/F)$ such that
$\tau|_E$ generates the cyclic group $\Gal(E/F)$.  Then
$\tau|_E$ has order $sp^m$, and there is $d\in\Z$ such that
$\tau(\zeta)=\zeta^d$.  It follows that the image of $d$ in
$(\Z/p^{m+1}\Z)^{\times}$ has order $sp^m$.  Since $\Gal(M/F)$
is abelian, $\tau(\alpha)/\alpha^d$ lies in the kernel of
$\N_{M/E}$.  Let $\rho$ be a generator for $\Gal(M/E)$.
Then by Hilbert's Theorem 90 there is $\beta\in M^{\times}$
such that $\tau(\alpha)/\alpha^d=\rho(\beta)/\beta$. 

     Let $\pi$ be a uniformizer for $M$, and write
$\beta=\gamma\pi^v$ with $\gamma\in\O_{M}^{\times}$ and
$v=v_M(\beta)$.  Set $\delta=\rho(\gamma)/\gamma$ and
$\epsilon=\rho(\pi)/\pi$, so that
$\tau(\alpha)/\alpha^d=\delta\epsilon^v$.
Now let $\alpha'$, $\gamma'$,
$\delta'$, $\epsilon'$ be elements of $\O_{M'}^{\times}$ which
correspond via $\nu$ to $\alpha$, $\gamma$, $\delta$,
$\epsilon$.  (In other words, we have $\nu(\alpha')\equiv\alpha
\pmod{\P_M^b}$, etc.)  In addition, choose $\pi'\in\P_{M'}$
such that $\theta(\pi')\equiv\pi\pmod{\P_M^{b+1}}$, and set
$\beta'=\gamma'{\pi'}^v$.  Let $\rho'=\jhat(\rho)$ be the element
of $\Gal(M'/K)$ which corresponds to $\rho\in\Gal(M/K)$.
Since $\rho(\pi)=\epsilon\pi$, it follows from (\ref{jhat})
that $\rho'(\pi')\equiv\epsilon'\pi'\pmod{\P_{M'}^{c+1}}$,
and hence that
\begin{equation}
\epsilon'\equiv\frac{\rho'(\pi')}{\pi'}\pmod{\P_{M'}^c}.
\end{equation}
Furthermore, since $\delta=\rho(\gamma)/\gamma$ and
$\delta\epsilon^v=\tau(\alpha)/\alpha^d$ we get
\begin{alignat}{2}
\delta'&\equiv\frac{\rho'(\gamma')}{\gamma'}&&\pmod{\P_{M'}^c} \\
\delta'{\epsilon'}^v&\equiv\frac{\tau'(\alpha')}{{\alpha'}^d}
&&\pmod{\P_{M'}^c}.
\end{alignat}
It follows that
\begin{equation} \label{modpb}
\frac{\tau'(\alpha')}{{\alpha'}^d}
\equiv\frac{\rho'(\gamma'{\pi'}^v)}{\gamma'{\pi'}^v}
\equiv\frac{\rho'(\beta')}{\beta'}
\pmod{\P_{M'}^{c}}.
\end{equation}

     Let $E'$ be the subfield of $M'$ fixed by
$\langle\rho'\rangle$.  Since $M=LE$, it follows from
Lemma~\ref{change} that the upper ramification breaks of
$M/E$ are bounded above by $\psi_{E/K}(u_{L/K})$.  Hence
the lower ramification breaks of $M/E$ are bounded above
by $\psi_{M/E}\circ\psi_{E/K}(u_{L/K})=\psi_{M/K}(u_{L/K})$,
which by assumption (\ref{cond3}) is less than $c$.  It
follows that the isomorphism between $\Gal(M/E)$ and $\Gal(M'/E')$
induced by $\jhat$ respects ramification filtrations, and
hence that $\psi_{M'/E'}=\psi_{M/E}$.  Thus by assumption
(\ref{cond2}) we have
$c>\psi_{M/E}(r)=\psi_{M'/E'}(r)$.  Therefore by (\ref{modpb})
and Proposition~\ref{norms} we get
\begin{equation} \label{tauzeta2}
\N_{M'/E'}\left(\frac{\tau'(\alpha')}{{\alpha'}^d}\right)
\equiv\N_{M'/E'}\left(\frac{\rho'(\beta')}{\beta'}\right)
\pmod{\P_{E'}^{r+1}}.
\end{equation}

     Let $\zeta'=\N_{M'/E'}(\alpha')$.  Since $\Gal(M'/K)$ is
abelian and $\rho'\in\Gal(M'/E')$, the congruence
(\ref{tauzeta2}) reduces to
\begin{equation} \label{zetaprime}
\frac{\tau'(\zeta')}{{\zeta'}^d}\equiv1\pmod{\P_{E'}^{r+1}}.
\end{equation}
(Note that if we simply defined $\zeta'$ to be an element of
$\O_E$ such that
$\nu(\zeta')\equiv\zeta\pmod{\P_M^b}$ then by (\ref{jhat})
and assumption (\ref{cond1}) we would get the weaker congruence
$\tau'(\zeta')
\equiv\zeta'^d\pmod{\P_{E'}^{q+1}}$.  This explains why we have
used such a roundabout method to define $\zeta'$.)
By applying Proposition~\ref{root} to (\ref{zetaprime}) with
$n=r+1>e_0(p^{m+1}+p^m-1)$ we get
$E'=F(\xi)$, with $\xi$ a primitive $p^{m+1}$th root of 1.
Therefore $E'=E$.
Furthermore, we have $v_E(\zeta'-\xi)\ge(gs+e_0)p^m$, where
\begin{equation}
g=\begin{cases}
\dst\left\lceil\frac{1+(y-e)s}{p^m}\right\rceil&
\text{if $h=0$ and $y>e$,} \\[.2cm]
1&\text{otherwise.}
\end{cases}
\end{equation}
Since $(gs+e_0)p^m>q$ we get $\xi\equiv\zeta'\pmod{\P_E^{q+1}}$.
By assumption (\ref{cond1}) we have $b\ge c>p^{n+t-m}q$, and hence
\begin{equation} \label{psieta3}
\nu(\xi)\equiv\nu(\zeta')\equiv
\nu(\N_{M'/E'}(\alpha'))\pmod{\P_M^{p^{n+t-m}q+1}}.
\end{equation}
Therefore by (\ref{jhat}) we get
\begin{equation} \label{psieta2}
\nu(\xi)\equiv\N_{M/E}(\alpha)\equiv\zeta
\pmod{\P_M^{p^{n+t-m}q+1}}.
\end{equation}

     Let $L_m/K$, $L_m'/K$ be the unique subextensions of
$L/K$, $L'/K$ of degree $p^m$, and set $M_m=L_mE=L_m(\zeta)$,
$M_m'=L_m'E=L_m'(\zeta)$.  Then $M_m/E$, $M_m'/E$ are the
unique subextensions of $M/E$, $M'/E$ of degree $p^t$.
Let $\pi_m=\N_{M/M_m}(\pi)$ and $\pi_m'=\N_{M'/M_m'}(\pi')$.
Then $\pi_m,\pi_m'$ are uniformizers for $M_m,M_m'$ such that
$\theta(\pi_m')\equiv\pi_m\pmod{\P_M^{c+1}}$.  Set
$\tilde{q}=\lfloor q/e_0\rfloor$ and
$c_m=e_0p^t\tilde{q}$.  By
assumption (\ref{cond1}) we have
\begin{equation}
c>p^{n+t-m}q=[M:M_m]p^tq\ge[M:M_m]c_m.
\end{equation}
Thus there is a unique $k$-linear ring homomorphism
\begin{equation}
\nu_m:\O_{M_m'}/\P_{M_m'}^{c_m}\lra\O_{M_m}/\P_{M_m}^{c_m}
\end{equation}
such that $\nu_m(\pi_m')\equiv\pi_m\pmod{\P_{M_m}^{c_m}}$
and a unique $\O_{M_m'}/\P_{M_m'}^{c_m}$-module homomorphism
\begin{equation}
\theta_m:\P_{M_m'}/\P_{M_m'}^{c_m+1}\lra\P_{M_m}/\P_{M_m}^{c_m+1}
\end{equation}
such that $\theta_m(\pi_m')\equiv\pi_m\pmod{\P_{M_m}^{c_m+1}}$.
These give an isomorphism $j_m=(1,\nu_m,\theta_m)$ from
$\Tr_{c_m}(M_m')$ to $\Tr_{c_m}(M_m)$.

     Let $\omega\in\Gal(\Q_p(\zeta)/\Q_p)$ be such that
$\omega(\xi)=\zeta$.
Since $\Z_p[\zeta]$ is the ring of integers of $\Q_p(\zeta)$,
it follows from (\ref{psieta2}) that
\begin{equation} \label{j0f}
j_m\circ f_{M_m'/\Q_p(\zeta)}
\equiv f_{M_m/\Q_p(\zeta)}\circ f_{\omega}\pmod{R(\tilde{q})}.
\end{equation}
Since $u_{L_m/K}=y+(m-h-1)e$ and $u_{K(\zeta)/K}=me$, we have
$u_{M_m/K}=y+(m-1)e$ if $h=0$ and $y>e$, and $u_{M_m/K}=me$
otherwise.  Therefore
\begin{align} \label{psiu}
\psi_{E/K}(u_{M_m/K})&=\begin{cases}
e_0(p^m-1)+sp^m(y-e)&\text{if $h=0$ and $y>e$}, \\
e_0(p^m-1)&\text{otherwise;}
\end{cases} \\
&=q-e_0. \label{psiuprime}
\end{align}

     Since $\psi_{M_m/K}(u_{M_m/K})$ and
$\psi_{M_m/E}(u_{M_m/E})$ are the largest lower ramification
breaks of $M_m/K$ and $M_m/E$, we have
$\psi_{M_m/K}(u_{M_m/K})\ge\psi_{M_m/E}(u_{M_m/E})$.
Applying $\phi_{M_m/E}$ to both sides of this inequality
we get $\psi_{E/K}(u_{M_m/K})\ge u_{M_m/E}$, and hence
$q-e_0\ge u_{M_m/E}$.  Since the ramification index $e_0$ of
$E/\Q_p(\zeta)$ is relatively prime to $p$, it follows from
Lemma~\ref{tame} that
\begin{equation}
u_{M_m/\Q_p(\zeta)}=e_0^{-1}\cdot u_{M_m/E}\le
e_0^{-1}(q-e_0)<\tilde{q}.
\end{equation}
Therefore $M_m/\Q_p(\zeta)$ and $M_m'/\Q_p(\zeta)$ satisfy
condition $C^{\tilde{q}}$.  Hence by applying
Corollary~\ref{induced} to (\ref{j0f}) we see
that $\omega$ can be extended to an isomorphism
$\tilde{\omega}:M_m'\ra M_m$ such that
\begin{equation} \label{j0f2}
j_m\equiv f_{\tilde{\omega}}\pmod{R(d)},
\end{equation}
where $d=\psi_{M_m/\Q_p(\zeta)}(\tilde{q})$.

     For $\gamma_m\in\Gal(M_m/K)$
let $\gamma$ be a lifting of $\gamma_m$ to $\Gal(M/K)$
and let $\jhat_m(\gamma_m)$ be the restriction of
$\jhat(\gamma)$ to $M_m'$.  Since
$\jhat(\Gal(M/M_m))=\Gal(M'/M_m')$ we see that
$\jhat_m(\gamma_m)$  does not depend on the choice of the
lifting $\gamma$.  Thus
\begin{equation}
\jhat_m:\Gal(M_m/K)\lra\Gal(M_m'/K)
\end{equation}
is a well-defined isomorphism.  By (\ref{jhat}) we have
\begin{equation} \label{jzero}
f_{\gamma_m}\circ j_m\equiv j_m\circ f_{\jhat_m(\gamma_m)}
\pmod{R(c_m)}.
\end{equation}
Since $c_m=e_0p^t\tilde{q}\ge d$, by (\ref{j0f2}) and
(\ref{jzero}) we get
\begin{alignat}{2} \label{omegatilde2}
f_{\gamma_m}\circ f_{\tilde{\omega}}&\equiv
f_{\tilde{\omega}}\circ f_{\jhat_m(\gamma_m)}&&\pmod{R(d)} \\
f_{\gamma_m}&\equiv f_{\tilde{\omega}\circ\jhat_m(\gamma_m)\circ\tilde{\omega}^{-1}}&&\pmod{R(d)}. \label{conj}
\end{alignat}

     Let $N/\Q_p$ be the smallest subextension of $M_m/\Q_p$
such that $M_m/N$ is Galois.  Then $\gamma_m$ and
$\tilde{\omega}\circ\jhat_m(\gamma_m)\circ\tilde{\omega}^{-1}$
both lie in $\Gal(M_m/N)$.  By (\ref{psiuprime}) we have
$\psi_{E/K}(u_{M_m/K})<e_0\tilde{q}=
\psi_{E/\Q_p(\zeta)}(\tilde{q})$.  It follows that the
largest lower ramification break $\psi_{M_m/K}(u_{M_m/K})$ of
$M_m/K$ is less than $\psi_{M_m/\Q_p(\zeta)}(\tilde{q})=d$.
Since $K/N$ is tamely ramified, by Lemma~\ref{tame} we see
that $\psi_{M_m/K}(u_{M_m/K})<d$ is also the largest lower
ramification break of $M_m/N$.  Hence by (\ref{conj}) we get
$\tilde{\omega}\circ\jhat_m(\gamma_m)\circ\tilde{\omega}^{-1}
=\gamma_m$ for all $\gamma_m\in\Gal(M_m/K)$.  Since
\begin{align}
\jhat_m(\Gal(M_m/K))&=\Gal(M_m'/K) \\
\jhat_m(\Gal(M_m/L_m))&=\Gal(M_m'/L_m')
\end{align}
this implies $\tilde{\omega}(K)=K$ and $\tilde{\omega}(L_m')=L_m$.
This proves Theorem~\ref{conditions}. \qed \medskip

     It remains to show that the values $a=ep^n$ and $m=m_0$
specified in Theorem~\ref{main} satisfy the inequalities
in Theorem~\ref{conditions}.  We prove this in the following
two lemmas.  The first of these lemmas, which is stronger
than needed to prove (\ref{cond1}), will also be used
in the proof of Theorem~\ref{proot}.

\begin{lemma} \label{cMEq}
Let $m\ge1$ satisfy $\psi_{L/K}((m+1+\frac{1}{p-1})e)<ep^n$.
Then $\psi_{M/L}(e(p^n-p^{n-1}))>p^{n+t-m}q$.
\end{lemma}

\proof Let $\tilde{a}=e(p^n-p^{n-1})$ and
$\tilde{c}=\psi_{M/L}(\tilde{a})$, and let
$\beta_{m-t},\beta_{m-t+1},\dots,\beta_{m-1}$ be the positive
upper ramification breaks of $M/L$.  Then we have
\begin{align}
\tilde{c}&=s\beta_{m-t}+sp(\beta_{m-t+1}-\beta_{m-t})+\dots
+sp^{t-1}(\beta_{m-1}-\beta_{m-2})+sp^t(\tilde{a}-\beta_{m-1}) \\[.2cm]
&=sp^t\tilde{a}-s(p-1)(\beta_{m-t}+p\beta_{m-t+1}+\dots
+p^{t-1}\beta_{m-1}). \label{tildec}
\end{align}
Since $E/\Q_p(\zeta)$ is tamely ramified, the positive
upper ramification breaks of $E/\Q_p$ are the same as the
positive upper ramification breaks $1,2,\dots,m$
of $\Q_p(\zeta)/\Q_p$.  It follows by Lemma~\ref{tame} that
the positive upper ramification breaks of $E/K$
are $(i+1)e$ for $0\le i<m$.  Therefore by
Lemma~\ref{change} we have $\beta_i\le\psi_{L/K}((i+1)e)$ for
$m-t\le i<m$.

     The values of $\psi_{L/K}((i+1)e)$ can be
computed using the ramification data for $L/K$ given in
(\ref{lowerLK}):
\begin{equation} \label{psiLKie}
\psi_{L/K}((i+1)e)=
\begin{cases}
\dst z+ep^{h+1}\cdot\frac{p^i-1}{p-1}+
p^{h+i+1}(e-y)&\text{if $y\le e$,} \\[.3cm]
\dst z+ep^{h+1}\cdot\frac{p^i-1}{p-1}+p^{h+i}(e-y)
&\text{if $y>e$.}
\end{cases}
\end{equation}
It follows from (\ref{tildec}) that
\begin{equation} \label{ccase1}
\tilde{c}\ge sp^t\tilde{a}+s(p^t-1)\left(\frac{ep^{h+1}}{p-1}-z\right)-
sp^{m+h-t+1}\cdot\frac{p^{2t}-1}{p+1}\left(\frac{p}{p-1}e-y\right)
\end{equation}
if $y\le e$, and
\begin{equation} \label{ccase2}
\tilde{c}\ge sp^m\tilde{a}+s(p^m-1)\left(\frac{ep^{h+1}}{p-1}-z\right)-
sp^h\cdot\frac{p^{2m}-1}{p+1}\left(\frac{2p-1}{p-1}e-y\right)
\end{equation}
if $y>e$.  In this last inequality we use the fact that $t=m$
if $y\not=e$.

     If $y>e$ then since $z\le p^hy$ and
$\dst p^m-1\le\frac{p^{2m}-1}{p+1}$, by (\ref{ccase2}) we get
\begin{align}
\tilde{c}&\ge sp^m\tilde{a}+s(p^m-1)\left(\frac{ep^{h+1}}{p-1}-p^hy\right)
-sp^h\cdot\frac{p^{2m}-1}{p+1}
\left(\frac{2p-1}{p-1}e-y\right) \\
&\ge sp^m\tilde{a}+s(p^m-1)\left(\frac{ep^{h+1}}{p-1}-p^he\right)
-sp^h\cdot\frac{p^{2m}-1}{p+1}\left(\frac{2p-1}{p-1}e-e
\right) \\
&=es\left(p^{n+m}-p^{n+m-1}+\frac{p^m-1}{p-1}\cdot
p^h-\frac{p^{2m}-1}{p^2-1}\cdot p^{h+1}\right). \label{ces}
\end{align}
By (\ref{mbound}) and the assumption 
$\psi_{L/K}((m+1+\frac{1}{p-1})e)<ep^n$
we have $m+h\le n-1$.  Therefore $\tilde{c}$ is greater than
\begin{align}
esp^{n+m}\left(1-\frac{1}{p}-\frac{1}{p^2-1}\right)
=\frac{p^3-p^2-2p+1}{p^2+p}\cdot e_0 p^{n+m}.
\end{align}
Since $q\le(\frac{es}{p-1}+e_0)p^m=2e_0p^m$ and $p\ge5$ we get
$\tilde{c}>2e_0p^{n+m}\ge p^nq$.  Since $t=m$ in this case
we conclude that $\tilde{c}>p^{n+t-m}q$.

     If $y\le e$ then since $h\le n-m-1\le n-2$ we have
$y>e/(p-1)$.  It follows by (\ref{ccase1}) that
\begin{align}
\tilde{c}&\ge sp^t\tilde{a}+s(p^t-1)\left(\frac{ep^{h+1}}{p-1}-p^hy\right)
-sp^{m+h-t+1}\cdot\frac{p^{2t}-1}{p+1}\left(\frac{p}{p-1}e-y\right) \\
&\ge sp^t\tilde{a}+s(p^t-1)\left(\frac{ep^{h+1}}{p-1}-
\frac{p^he}{p-1}\right)
-sp^{m+h-t+1}\cdot\frac{p^{2t}-1}{p+1}\left(\frac{p}{p-1}e-\frac{e}{p-1}\right) \\
&=es\left(p^{n+t}-p^{n+t-1}+(p^t-1)p^h-\frac{p^{2t}-1}{p+1}\cdot p^{m+h-t+1}\right).
\end{align}
As above this implies that $\tilde{c}$ is greater than or
equal to
\begin{equation}
esp^{n+t}\left(1-\frac{1}{p}-\frac{1}{p+1}\right)
=\frac{p^3-2p^2+1}{p^2+p}\cdot e_0p^{n+t},
\end{equation}
and hence that $\tilde{c}>e_0p^{n+t}= p^{n+t-m}q$. \qed \medskip

     It follows from Lemma~\ref{cMEq} that assumption (\ref{cond1})
is satisfied by $a=ep^n$, $m=m_0$.  We now show that these
values satisfy assumptions (\ref{cond2}) and (\ref{cond3}) as
well.

\begin{lemma} \label{psid}
Let $a,m\ge1$ satisfy $\frac{2}{p-1}\cdot ep^n<a\le
ep^n$ and $\psi_{L/K}((m+1+\frac{1}{p-1})e)<a$.  Then
$\psi_{M/L}(a)>\psi_{M/E}(r)$ and $\psi_{M/L}(a)>\psi_{M/K}(u_{L/K})$.
\end{lemma}

\proof Since the positive upper ramification breaks of $E/K$
are $(i+1)e$ for $0\le i<m$, the positive lower ramification
breaks
of $E/K$ are $\psi_{E/K}((i+1)e)=e_0(p^{i+1}-1)$ for $0\le i<m$.
It follows that
\begin{equation} \label{psiEKd}
\phi_{E/K}(r)=\begin{cases}
\dst \left(m+1+\frac{1}{p-1}\right)e+(y-e)&
\text{if $h=0$ and $y>e$}, \\[.4cm]
\dst \left(m+1+\frac{1}{p-1}\right)e&\text{otherwise}.
\end{cases}
\end{equation}
If $h=0$ and $y>e$ then
\begin{align} \label{phiepn}
\phi_{L/K}(a)=y+(n-1)e+\frac{1}{p^n}\left(a-\left(y+ep\cdot
\frac{p^{n-1}-1}{p-1}\right)\right)
\end{align}
is greater than $\phi_{E/K}(r)$, since $m\le n-1$ and
$a>\frac{2}{p-1}\cdot ep^n$.  In the other
cases the inequality ${\phi_{L/K}(a)>\phi_{E/K}(r)}$
follows from the assumption
$\psi_{L/K}((m+1+\frac{1}{p-1})e)<a$ and (\ref{psiEKd}).
Applying $\psi_{M/K}$ to both sides of this
inequality we get $\psi_{M/L}(a)>\psi_{M/E}(r)$.  By
(\ref{lowerLK}) we have
\begin{equation}
\psi_{L/K}(u_{L/K})=z+ep^{h+1}\cdot\frac{p^{n-h-1}-1}{p-1}.
\end{equation}
Since $z\le p^hy\le ep^{h+1}/(p-1)$, this quantity
is less than $a$.  It follows that
$\psi_{M/K}(u_{L/K})<\psi_{M/L}(a)$. \qed \medskip

     Theorem~\ref{main} follows from
Theorem~\ref{conditions} combined with Lemmas~\ref{cMEq} and
\ref{psid}.  To prove Theorem~\ref{proot} we apply
Theorem~\ref{conditions} to the subextensions
of $L/K$ and $L'/K$ of degree $p^{n-1}$:\medskip

\noindent {\em Proof of Theorem~\ref{proot}:}
Let $L/K$, $L'/K$ be totally ramified $(\Z/p^n\Z)$-extensions
which satisfy condition (*) of Theorem~\ref{main}.  We may
assume without loss of generality that $K$ contains a
primitive $p$th root of unity, that $m_0\ge2$, and that
$n\ge3$.  Since $p\nmid e$ we see that $K$ contains no primitive
$p^2$th roots of unity.  Therefore the group $\mu$ of
$p$-power roots of unity of $K$ is cyclic of order $p$.
Let $L_{n-1}/K$ be the unique
subextension of $L/K$ of degree $p^{n-1}$.  Then
$\N_{L/K}(L^{\times})$ has index $p$ in
$\N_{L_{n-1}/K}(L_{n-1}^{\times})$, so
$\mu\le\N_{L_{n-1}/K}(L_{n-1}^{\times})$.
Hence by Lemma~\ref{Zp} we see that $L_{n-1}$
is contained in a $\Z_p$-extension $L_{\infty}$ of $K$.
Let $j=\psi_{L/K}(u_{L/K})$ be the unique ramification break
of $L/L_{n-1}$ and let $l=\lceil\frac{p-1}{p}\cdot j\rceil$.
Then by (\ref{lowerLK}) we have
\begin{equation} \label{lform}
l=\left\lceil\frac{p-1}{p}\cdot\left(z+ep^{h+1}\cdot
\frac{p^{n-h-1}-1}{p-1}\right)\right\rceil.
\end{equation}
It follows from \cite[Prop.\:2.2.1]{cn} that the norm map
induces ring isomorphisms
\begin{align}
\overline{\N}_{L/L_{n-1}}&:\O_L/(\P_L^l)\lra
\O_{L_{n-1}}/(\P_{L_{n-1}}^l) \label{Nbar} \\
\overline{\N}_{L'/L_{n-1}'}&:\O_L/(\P_{L'}^l)\lra
\O_{L_{n-1}'}/(\P_{L_{n-1}'}^l). \label{Nbarprime}
\end{align}
These isomorphisms are Galois-equivariant and induce the
$p$-Frobenius map on $k$.

     Let $\sigma_{n-1}$ denote the restriction of $\sigma$ to
$L_{n-1}$, set $\pi_{L_{n-1}}=\N_{L/L_{n-1}}(\pi_L)$, and
set $\pi_{L_{n-1}'}=\N_{L'/L_{n-1}'}(\pi_{L'})$.  By applying
the arguments used to prove (\ref{normcong}) to
(\ref{Nbar}) and (\ref{Nbarprime}) we get
\begin{alignat}{2}
h_{\pi_{L_{n-1}}}^{\sigma_{n-1}}(X)&\equiv
(h_{\pi_L}^{\sigma})^{\phi}(X)&&\pmod{X^l} \\
h_{\pi_{L_{n-1}'}}^{\sigma_{n-1}'}(X)&\equiv
(h_{\pi_{L'}}^{\sigma'})^{\phi}(X)&&\pmod{X^l}.
\end{alignat}
Since $h_{\pi_L}^{\sigma}=h_{\pi_{L'}}^{\sigma'}$ this
implies
\begin{equation}
h_{\pi_{L_{n-1}}}^{\sigma_{n-1}}(X)\equiv
h_{\pi_{L_{n-1}'}}^{\sigma_{n-1}'}(X)\pmod{X^l}.
\end{equation}

     Since $m_0\ge2$ we have $h\le n-3$, so $y$ is the smallest
upper ramification break of $L_{n-1}/K$ which exceeds
$\frac{1}{p-1}\cdot e$.  Therefore the largest upper
ramification break of $L_{n-1}/K$ is $u_{L_{n-1}/K}=y+(n-h-2)e$.
Let $m=m_0-1$ and define $E/F$ as in the proof of
Theorem~\ref{conditions}.  Also set $M_{n-1}=L_{n-1}E$.
To prove
Theorem~\ref{proot} it suffices by Theorem~\ref{conditions}
to prove the following inequalities:
\begin{align}
\psi_{M_{n-1}/L_{n-1}}(l)&>[M_{n-1}:E]q \label{newcond1} \\[.2cm]
\psi_{M_{n-1}/L_{n-1}}(l)&>\psi_{M_{n-1}/E}(r) \label{newcond2} \\[.2cm]
\psi_{M_{n-1}/L_{n-1}}(l)&>\psi_{M_{n-1}/K}(u_{L_{n-1}/K}). \label{newcond3}
\end{align}

     Using (\ref{lowerLK}) to compute the left side of the
inequality $\psi_{L/K}((m_0+1+\frac{1}{p-1})e)<ep^n$ we get
\begin{equation} \label{epn}
z+ep^{h+1}\cdot\frac{p^{m_0}-1}{p-1}+p^{m_0+h+1}\left(
\frac{p}{p-1}\cdot e-y\right)<ep^n.
\end{equation}
Since $z\le p^hy\le ep^{h+1}/(p-1)$ we have $(p-1)z-ep^{h+1}\le0$.
Adding this inequality to (\ref{epn}) and dividing by $p$ gives
\begin{equation}
z+ep^{h+1}\cdot\frac{p^{m_0-1}-1}{p-1}+p^{m_0+h}\left(
\frac{p}{p-1}\cdot e-y\right)<ep^{n-1}.
\end{equation}
Hence we have 
\begin{equation}
\psi_{L_{n-1}/K}\left(\left(m+1+\frac{1}{p-1}\right)e\right)=
\psi_{L_{n-1}/K}\left(\left(m_0+\frac{1}{p-1}\right)e\right)<ep^{n-1}.
\end{equation}
It follows from (\ref{lform}) that
$l>e(p^{n-1}-p^{n-2})$.  Therefore
(\ref{newcond1}) follows from Lemma~\ref{cMEq}.  By
(\ref{lowerLK}) we have
\begin{equation}
\psi_{L_{n-1}/K}(u_{L_{n-1}/K})=z+ep^{h+1}\cdot\frac{p^{n-h-2}-1}{p-1}.
\end{equation}
Using (\ref{lform}) and the inequality $z\le ep^{h+1}/(p-1)$ we
deduce that $\psi_{L_{n-1}/K}(u_{L_{n-1}/K})<l$.  Applying
$\psi_{M_{n-1}/L_{n-1}}$ to this last inequality gives
(\ref{newcond3}).

     It remains to prove (\ref{newcond2}).  If $h=0$ and
$y>e$ then by (\ref{lform}) we have $l>(p^{n-1}-1)e$.  It
follows using (\ref{lowerLK}) that
\begin{equation}
\phi_{L_{n-1}/K}(l)>\left(n-1-\frac{1}{p-1}+
\frac{1}{p^n-p^{n-1}}\right)e+
\left(1-\frac{1}{p^{n-1}}\right)y.
\end{equation}
By (\ref{psiEKd}) with $m=m_0-1$ we have
\begin{equation}
\phi_{E/K}(r)=\left(m_0+\frac{1}{p-1}\right)e+(y-e).
\end{equation}
Since $m_0\le n-1$, $y\le(1+\frac{1}{p-1})e$, $p\ge5$, and
$n\ge2$, we get
$\phi_{E/K}(r)<\phi_{L_{n-1}/K}(l)$.  Applying
$\psi_{M_{n-1}/K}$ to this inequality gives
(\ref{newcond2}) in this case.

     Suppose $h\ge 1$ or $y\le e$.  Adding
\begin{equation}
(p-1)z-ep^{h+1}\le p\left\lceil\frac{p-1}{p}\cdot z\right\rceil-ep^{h+1}.
\end{equation}
to (\ref{epn}) and dividing by $p$ gives
\begin{equation}
z+ep^{h+1}\cdot\frac{p^{m_0-1}-1}{p-1}+p^{m_0+h}\left(
\frac{p}{p-1}\cdot e-y\right)<\left\lceil\frac{p-1}{p}\cdot
z\right\rceil+ep^{n-1}-ep^h.
\end{equation}
It follows from (\ref{psiEKd}), (\ref{lowerLK}), and
(\ref{lform}) that this inequality can be rewritten as
$\psi_{L_{n-1}/K}\circ\phi_{E/K}(r)<l$.
By applying $\psi_{M_{n-1}/L_{n-1}}$ we get
(\ref{newcond2}). \qed

\section{$p$-adic dynamical systems} \label{dynam}

     Let $K$ be a finite extension of $\Q_p$ and let $\P^{alg}$
be the maximal ideal in the ring of integers of $\Q_p^{alg}$.
Let $u(X)\in\O_K[[X]]$ be a power series such that
$u(0)=0$ and $u'(0)$ is a 1-unit.
We are interested in studying the periodic points of
$u(X)$.  These are the elements $\alpha\in\P^{alg}$
such that $u^{\circ m}(\alpha)=\alpha$ for some $m\ge1$;
the smallest such $m$ is called the period of $\alpha$.
Since $u'(0)$ is a 1-unit, it follows from
\cite[Cor.\:2.3.2]{pper} that all periodic points of $u(X)$
have period $p^n$ for some $n\ge0$.
In the introduction to \cite{nonarch}, Lubin stated
that the extension fields generated by the periodic points of
$u(X)$ are ``almost completely unknown''.
In this section we show how Theorem~\ref{conditions}
can be used to study the extension $K(\alpha)/K$ generated by
a single periodic point $\alpha$.

     Let $\ubar(X)\in k[[X]]$ denote the reduction of
$u(X)$ modulo $\P_K$.  It follows from our assumptions that
$\ubar(X)$ is an element of the group $\A(k)$ which was
defined in Section~\ref{norm}.  For $n\ge0$ let $i_n$ denote the
depth of $\ubar^{\circ p^n}(X)$.  If $i_n<\infty$ then
$i_n+1$ is equal to the number of solutions in $\P$ to the
equation $u^{\circ p^n}(X)=X$, counted with multiplicity.  
Let $\Gamma$ be the closed subgroup of $\A(k)$ generated by
$\ubar(X)$ and assume that $\Gamma$ is infinite; then
$\Gamma\cong\Z_p$.  It follows from Proposition~\ref{bij}
that there is a local field $L_0$ with residue field $k$,
a totally ramified $\Z_p$-extension $L_{\infty}/L_0$, and a
compatible
sequence of uniformizers $(\pi_n)$ for $L_{\infty}/L_0$
such that $\Gamma=\Gamma_{L_{\infty}/L_0}^{(\pi_n)}$.
The extension $L_{\infty}/L_0$
is determined uniquely up to $k$-isomorphism by $\Gamma$.
By \cite[Cor.\:3.3.4]{cn} the ramification data of the
extension $L_{\infty}/L_0$ is the same as the ramification
data of $\Gamma$.  We define the index $d$ of $\Gamma$ to be
the absolute ramification index of $L_0$; if $L_0$ has
characteristic $p$ then the index of $\Gamma$ is $\infty$.
If $d<\infty$ then it follows from Lemma~\ref{breaks}
that $b_n-b_{n-1}=d$ for all sufficiently large $n$.

\begin{theorem} \label{ext}
Let $p>3$, let $1\le d\le p-2$, and let $K/\Q_p$ be a finite
extension with
ramification index $e\le p-1$.  Then there is a finite tamely
ramified extension $E/K$ with the following property: Let
$u(X)\in\O_K[[X]]$ be a power series such that the closed
subgroup $\Gamma$ of $\A(k)$ generated by $\ubar(X)$
is isomorphic to $\Z_p$ and has index $d$.  Let $L_0$
be a local field with residue field $k$ and let
$L_{\infty}/L_0$ be a totally ramified $\Z_p$-extension
such that $\Gamma\in[\Gamma_{L_{\infty}/L_0}]$.
For $n\ge1$ let $L_n/L_0$ denote the subextension of
$L_{\infty}/L_0$ of degree $p^n$.
Then for each periodic point $\alpha$ of $u(X)$ with period
$p^n$ there is an embedding $\omega:L_n\ra\Q_p^{alg}$ such that 
\begin{equation}
[E(\alpha)\cap(E\cdot\omega(L_n)):E]\ge p^{n-2}.
\end{equation}
\end{theorem}

     Thus when the hypotheses of Theorem~\ref{ext} are satisfied
the special fiber $\ubar(X)$ of $u(X)$ carries a large amount
of information about the field extensions generated by
periodic points of $u(X)$.  It follows from Lemma~\ref{breaks}
that if the index of $\Gamma$ is $\infty$ then the upper
ramification breaks of $\Gamma$ satisfy $b_n\ge pb_{n-1}$ for
all $n\ge1$, while if the index of $\Gamma$ is $d<\infty$ then
for each $n\ge1$ we have either $b_n\ge pb_{n-1}$ or
$b_n-b_{n-1}=d$.  Therefore the index of $\Gamma$
can be effectively computed as long as it is finite.

     The rest of this section is devoted to proving
Theorem~\ref{ext}.  Let $l$ be the field extension of $k$ of
degree $d!$ and let $F$ be the unramified extension
of $\Q_p$ with residue field $l$.  The field $E$ is defined
to be the compositum of all totally ramified extensions $R/F$
of degree $d!e$.  The following lemma is a consequence of the
well-known properties of tamely ramified extensions of a
local field.

\begin{lemma} \label{facts}
(a) The absolute ramification index of $E$ is $d!e$. \\
(b) If $M$ is a finite extension of $\Q_p$ whose absolute
ramification index divides $d!e$ and whose residue field is
contained in $l$, then $M$ is contained in $E$.
\end{lemma}

     It follows from this lemma that $E$ is an extension of
$K$ with ramification index $d!$, and hence that $E/K$ is
tamely ramified.  In particular, if $d=1$ then $E/K$ is
unramified.

     Before proving Theorem~\ref{ext} we study the basic
properties of periodic points of power series which
satisfy the hypotheses of the theorem.  In particular, we are
interested in the degrees and ramification indices of extensions
generated by these periodic points.  For the remainder of
this section we assume without loss of generality that
$n\ge3$.

     As above we let $(i_n)_{n\ge0}$ and $(b_n)_{n\ge0}$
denote the lower and upper ramification sequences of $\Gamma$.
By Lemma~\ref{breaks}(a)
we have $b_0\le(1+\frac{1}{p-1})d$.  Since $d<p-1$ this implies
$i_0=b_0\le d$.  We also have $b_0\ge1>d/(p-1)$.  Therefore by
Lemma~\ref{breaks}(c) we get $b_n=b_{n-1}+d$ and hence
$i_n=i_{n-1}+dp^n$ for
all $n\ge1$.  We may write $u(X)=a_0X+a_1X^2+a_2X^3+\cdots$
with $a_i\in\O_K$ and $a_0$ a 1-unit.  Since $p>e/(p-1)$ we have
\begin{equation} \label{vK}
v_K(a_0^{p^n}-1)=v_K(a_0^{p^{n-1}}-1)+e
\end{equation}
if $a_0^p\not=1$.  It follows from \cite[Cor.\:2.3.1]{pper}
that $u^{\circ p^{n-1}}(X)-X$ divides $u^{\circ p^n}(X)-X$ in
$\O_K[[X]]$.  Let
\begin{equation}
q_n(X)=\frac{u^{\circ p^n}(X)-X}{u^{\circ p^{n-1}}(X)-X}.
\end{equation}
If $a_0^p\not=1$ then by
(\ref{vK}) the constant term of $q_n(X)$ has
$p$-valuation 1, while if $a_0^p=1$ then
the constant term of $q_n(X)$ is equal to $p$.  Note that
the Weierstrass degree of $q_n(X)\in\O_K[[X]]$ is
$i_n-i_{n-1}=dp^n$.

     Let $\alpha\in\P$ be a periodic point of $u(X)$ with period
$p^n$ and let
$M/K$ be the Galois closure of $K(\alpha)/K$.  Let
$G=\Gal(M/K)$, let
\begin{equation}
H=\{\tau\in G:\tau(\alpha)=u^{\circ i}(\alpha)
\text{ for some }i\in\Z\},
\end{equation}
and let $\sigma_1,\dots,\sigma_h$ be coset representatives
for $G/H$.  The polynomial
\begin{equation}
f(X)=\prod_{j=1}^h\prod_{i=0}^{p^n-1}
(X-u^{\circ i}(\sigma_j(\alpha)))
\end{equation}
lies in $\O_K[X]$ and has distinct roots, all of which
are zeros of $q_n(X)$.  Therefore $f(X)$ divides
$q_n(X)$ in $\O_K[[X]]$, and hence the
constant term $c$ of $f(X)$ is an element of $\P_K$ which
divides $p$.  It follows that $c$ has $p$-valuation $s/e$
for some $1\le s\le e$.  Since each of the $hp^n$ roots of
$f(X)$ has the same $p$-valuation as $\alpha$, we get
$v_p(\alpha)=s/ehp^n$.

     For each $1\le j\le h$ the set
$B_j=\{u^{\circ i}(\sigma_j(\alpha)):0\le i<p^n\}$
is a block for the permutation representation of $G$ acting
on the roots of $f(X)$.  Let $N$ be the kernel of the action
of $G$ on the set of blocks, let $T$ be the fixed field of
$N$, and let $U/K$ be the maximum unramified subextension of
$T/K$.  Since the degree of $f(X)$ is less than or equal to the
Weierstrass degree of $q_n(X)$ we have $h\le d$.  Since
$\Gal(T/K)\cong G/N$ is isomorphic to a
subgroup of $S_h$ this implies that $[T:U]$ and $[U:K]$
both divide $d!$.  Therefore $T$ is an extension of $\Q_p$
whose absolute ramification
index divides $d!e$ and whose residue field is contained in
$l$.  Hence by Lemma~\ref{facts}(b), $T$ is contained in $E$.

     For each $\tau\in N=\Gal(M/T)$ there is a unique
$i\in\Z/p^n\Z$ such that $\tau(\alpha)=u^{\circ i}(\alpha)$.
Hence $T(\alpha)/T$ is Galois, and $\Gal(T(\alpha)/T)$
can be identified with a subgroup of
$\Z/p^n\Z$.  It follows that $E(\alpha)/E$ is also Galois,
with $\Gal(E(\alpha)/E)$ isomorphic to a subgroup of
$\Gal(T(\alpha)/T)$.  Since the ramification
index of $E/\Q_p$ is $d!e$, the $E$-valuation of $\alpha$
is $t/p^n$, where $t=(d!/h)\cdot s$ is relatively prime to
$p$.  It follows that $[E(\alpha):E]\ge p^n$, and hence that
$\Gal(E(\alpha)/E)\cong\Z/p^n\Z$.

     Since the absolute ramification index of $T$ divides
$d!e$ and the residue field of $T$ is contained in $l$, there
is a totally ramified extension $R/F$ of degree $d!e$ such
that $R$ contains $T$.  Then $E(\alpha)$ is an unramified
extension of $R(\alpha)$, so the
$R(\alpha)$-valuation of $\alpha$ is $t$.  Therefore we
can write $\alpha=\zeta\pi^t$, where $\zeta\in F$ is a root
of unity whose order is prime to $p$ and $\pi$ is a
uniformizer for $R(\alpha)$.  Let $K(l)=FK$ be the unramified
extension of $K$ with residue field $l$ and let
$\tau\in\Gal(E(\alpha)/E)$ satisfy $\tau(\alpha)=u(\alpha)$.
Let $v(X)$ be the unique element of 
$\O_{K(l)}[[X]]$ such that $\zeta v(X)^t=u(\zeta X^t)$ and
$v'(0)\equiv1\pmod{\P_{K(l)}}$.  Then
\begin{equation}
\zeta
v(\pi)^t=u(\zeta\pi^t)=u(\alpha)=\tau(\alpha)=\zeta\tau(\pi)^t.
\end{equation}
Since $\tau$ has order $p^n$ this implies $\tau(\pi)=v(\pi)$.
We are now in a position to prove the following key fact:

\begin{prop} \label{td}
The absolute ramification index $d!e$ of $E$ is equal to $td$.
\end{prop}

\proof Let $\Gamma'\cong\Z_p$ be
the closed subgroup of $\A(l)$ generated by $\vbar(X)$.
Since $E(\alpha)=E(\pi)$ and $\tau(\pi)=v(\pi)$,
the lower ramification breaks of the extension
$E(\alpha)/E$ which are less than $d!p^n$ (the ramification
index of $E(\alpha)/K(l)$) are the same as the
lower ramification breaks of $\Gamma'$ which are less than
$d!p^n$.  It follows from the definition of $v(X)$ that the
ramification breaks of $\Gamma'$ are $t$ times
the ramification breaks of $\Gamma$.  Therefore the first
two lower ramification breaks of $\Gamma'$ are $ti_0$ and
$ti_1=ti_0+tdp$.  Using the inequalities $s\le p-1$ and
$i_0\le d\le p-2$ we deduce that $ti_0+tdp<d!p^3$.
Therefore the first two lower ramification breaks of
$E(\alpha)/E$ are $ti_0$ and $ti_0+tdp$, and hence the
first two upper ramification breaks of $E(\alpha)/E$ are
$ti_0$ and $ti_0+td$.  Since $ti_0+td<p\cdot ti_0$, it
follows from Lemma~\ref{breaks} that the absolute ramification
index of $E$ is $td$. \qed

\begin{cor} \label{same}
For $n\ge3$ the periodic points of $u(X)$ with period $p^n$
all have $p$-valuation $1/dp^n$.
\end{cor}

\proof Let $\alpha$ be a periodic point of $u(X)$ with period
$p^n$.  We saw above that $v_p(\alpha)=s/ehp^n$.  Since
$s=ht/d!$ and $t=d!e/d$ we get $v_p(\alpha)=1/dp^n$. \qed

\begin{prop} \label{fact}
Let $n\ge3$.  Then every zero of $q_n(X)$ is
periodic with period $p^n$.
\end{prop}

\proof Let $\alpha$ be a periodic point with period $p^j$ for
some $0\le j<n$.  If $j=0$ then $\alpha$ is a zero of $u(X)-X$,
and if $j\ge1$ then $\alpha$ is a zero of $q_j(X)$.  It
follows by the Weierstrass preparation theorem that $\alpha$
is a root of a distinguished polynomial
with coefficients in $\O_K$ which divides $u(X)$ or $q_j(X)$.
Since $u(X)$ has Weierstrass degree $i_0+1$, and $q_j(X)$ has
Weierstrass degree $dp^j$, we must have
\begin{equation}
v_p(\alpha)\ge
\begin{cases}
\dst\frac{1}{e(i_0+1)}&\text{if $j=0$,} \\[.4cm]
\dst\frac{1}{edp^j}&\text{ if $1\le j<n$.}
\end{cases}
\end{equation}
In particular, since $n\ge3$, $e<p$, and $i_0\le d$ we have
$v_p(\alpha)>1/dp^n$.

     The series $q_n(X)$ has $dp^n$ zeros, counting
multiplicities; all of these are periodic points with period
$p^j$ for some $0\le j\le n$.  Let $\alpha\in\P^{alg}$ be a
zero of $q_n(X)$.  If $\alpha$ is a periodic point
with period $p^n$ then by Corollary~\ref{same} we have
$v_p(\alpha)=1/dp^n$.  On the other hand, if $\alpha$ is a
periodic point with period $p^j$ for some $0\le j<n$, then
$v_p(\alpha)>1/dp^n$.  The sum of the $p$-valuations
of the $dp^n$ zeros of $q_n(X)$ is 1.  Therefore
all the zeros of $q_n(X)$ must have period $p^n$.
\qed \medskip

     It follows that for $n\ge3$ the
periodic points of $u(X)$ with period $p^n$ are precisely
the zeros of $q_n(X)$, and that the number of periodic
points of $u(X)$ of period $p^n$, counted with multiplicity,
is equal to the Weierstrass degree $i_n-i_{n-1}=dp^n$ of
$q_n(X)$.  In particular, $u(X)$ has periodic points of
period $p^n$ for every $n\ge3$. \medskip

\noindent {\em Proof of Theorem~\ref{ext}:} Since
$\Gamma\in[\Gamma_{L_{\infty}/L_0}]$, there exists
a compatible sequence of uniformizers
$(\pi_j)$ for $L_{\infty}/L_0$
such that $\Gamma=\Gamma_{L_{\infty}/L_0}^{(\pi_j)}$.  Since
$\ubar(X)$ generates $\Gamma$, it follows from (\ref{gsig})
that there is
a generator $\sigma$ for $\Gal(L_{\infty}/L_0)$ such that
\begin{equation} \label{sigpi}
\sigma(\pi_j)\equiv\ubar^{\phi^{-j}}(\pi_j)\pmod{\P_{L_j}^{r_j+1}}
\end{equation}
for all $j\ge1$,
where $r_j=\lceil(p-1)i_j/p\rceil$, and we identify $k$
with a subring of $\O_{L_j}/(\pi_j^{r_j+1})$ using the
Teichm\"uller lifting.

     The map $x\mapsto x^p$ is an automorphism of the group
of roots of unity of $F$.  We denote the inverse of this
automorphism by raising to the power $p^{-1}$.  For $1\le
j\le\infty$ let $E_j=EL_j$.

\begin{lemma} \label{seq}
There exists a compatible sequence of uniformizers
$(\tilde{\pi}_j)$ for $E_{\infty}/E$ such that
$\pi_j=\zeta^{p^{-j}}{\tilde{\pi}_j}^t$ for $0\le j<\infty$.
\end{lemma}

\proof Let $j\ge0$ and let $\tilde{\pi}_j\in\Q_p^{alg}$ be a
root of $X^t-\zeta^{-p^{-j}}\pi_j$.  Let $k_E$ denote
the residue field of $E$ and let $E_j'$ be the unramified
extension of $FL_j(\tilde{\pi}_j)$ with residue field $k_E$.
Since $\zeta^{-p^{-j}}\pi_j$ is a uniformizer for $FL_j$, the
extension of $FL_j(\tilde{\pi}_j)/FL_j$ is totally ramified,
with ramification index
$t$.  Therefore the maximum tame subextension $T_j/\Q_p$ of
$FL_j(\tilde{\pi}_j)/\Q_p$ has ramification index $td=d!e$
and residue field $l$.  It follows by Lemma~\ref{facts}(b)
that $T_j$ is contained in $E$.  Thus
by Lemma~\ref{facts}(a), $E$ is an unramified
extension of $T_j$, so $E$ is contained in $E_j'$.  Since
$L_j\subset E_j'$, we get $E_j'=EL_j=E_j$.  The norm
map $\N_{E_j/E}$ gives a bijection between the roots of
$X^t-\zeta^{-p^{-j}}\pi_j$ and the roots of
$X^t-\zeta^{-1}\pi_0$.  Therefore we may assume that
$\N_{E_j/E}(\tilde{\pi}_j)=\tilde{\pi}_0$ for every $j\ge1$.  It
follows from this assumption that $(\tilde{\pi}_j)_{j\ge0}$
is a compatible
sequence of uniformizers for $E_{\infty}/E$. \qed \medskip

     Since $\zeta v(X)^t=u(\zeta X^t)$ we have
$\zetabar\vbar(X)^t=\ubar(\zetabar X^t)$, where $\zetabar$
denotes the image of $\zeta$ in $l\cong\O_{K(l)}/\P_{K(l)}$.
Applying $\phi^{-n}$ we get
$\zetabar^{p^{-n}}\vbar^{\phi^{-n}}(X)^t=
\ubar^{\phi^{-n}}(\zetabar^{p^{-n}}X^t)$.
Let $\tilde{\sigma}$ be the generator for
$\Gal(E_{\infty}/E)$ whose restriction to
$L_{\infty}$ is $\sigma$.  Then by (\ref{sigpi})
and Lemma~\ref{seq} we have
\begin{alignat}{2}
\tilde{\sigma}(\zeta^{p^{-n}}{\tilde{\pi}_n}^t)&\equiv
\ubar^{\phi^{-n}}(\zeta^{p^{-n}}{\tilde{\pi}_n}^t)&&
\pmod{\P_{E_n}^{t(r_n+1)}} \\
\zeta^{p^{-n}}\tilde{\sigma}({\tilde{\pi}_n})^t&\equiv
\zeta^{p^{-n}}\vbar^{\phi^{-n}}({\tilde{\pi}_n})^t&&
\pmod{\P_{E_n}^{t(r_n+1)}} \\
\tilde{\sigma}(\tilde{\pi}_n)&\equiv\vbar^{\phi^{-n}}(\tilde{\pi}_n)&&
\pmod{\P_{E_n}^{tr_n+1}}.
\label{twisted}
\end{alignat}
Let $\Phi$ be an automorphism of $\Q_p^{alg}$ which
induces the $p$-Frobenius on residue fields, and let
$\Theta:E_{\infty}\ra\Phi^n(E_{\infty})$
be the isomorphism induced by $\Phi^n$.  Applying
$\Theta$ to (\ref{twisted}) we get
\begin{equation} \label{untwisted}
\hat{\sigma}(\hat{\pi}_n)\equiv\vbar(\hat{\pi}_n)
\pmod{\P_{\Theta(E_n)}^{tr_n+1}},
\end{equation}
where $\hat{\pi}_n=\Theta(\tilde{\pi}_n)$ is a uniformizer for
$\Theta(E_n)$ and
$\hat{\sigma}=\Theta\circ\tilde{\sigma}\circ\Theta^{-1}$
is a generator for $\Gal(\Theta(E_{\infty})/E)$.
(Note that since $E$ is Galois over $\Q_p$ we
have $\Theta(E)=E$.)
On the other hand, since $\tau(\pi)=v(\pi)$ we have
\begin{equation} \label{taucong}
\tau(\pi)\equiv\vbar(\pi)\pmod{\P_E^{tr_n+1}}.
\end{equation}
Note that since $\pi$ is a uniformizer for $R(\alpha)$, $\pi$
is also a uniformizer for $E$.  To complete the proof of
Theorem~\ref{ext} we will
apply Theorem~\ref{conditions} to the extensions
$\Theta(E_n)/E$ and $E(\alpha)/E$.  To do this we must first
compute some ramification data.

     Since $d<p-1$, it follows from Lemma~\ref{breaks}(b)
that the $j$th upper ramification break of $L_{\infty}/L_0$
is $b_j=b_0+jd$.  Therefore the lower breaks of
$L_{\infty}/L_0$ are given by $i_j=i_0+dp+dp^2+\dots+dp^j$,
with $i_0=b_0$.
The unique ramification break of $L_{n+1}/L_n$ is equal to
the $n$th lower ramification break $i_n$ of $L_{\infty}/L_0$.
It follows that
\begin{align}
r_n&=\left\lceil\frac{p-1}{p}\cdot(i_0+dp+dp^2+\dots+dp^n)
\right\rceil \\
&>d(p^n-1). \label{snbound}
\end{align}

     The ramification breaks of $E_n/E$ are $t$ times
the ramification breaks of $L_n/L_0$.  The upper and lower
ramification breaks of $L_n/L_0$ are the integers
$b_j$ and $i_j$ for $0\le j<n$ which were computed in the
preceding paragraph.  Therefore we have
\begin{align}
\psi_{E_n/E}\left(\left(n-1+\frac{1}{p-1}\right)td\right)&=
\frac{p^n+p^{n-1}-p}{p-1}\cdot td-(p^{n-1}-1)ti_0.
\end{align}
This value is less than $td(p^n-p^{n-1})$, which by
(\ref{snbound}) is less than $tr_n$.
Comparing (\ref{untwisted}) with (\ref{taucong}) and
applying Lemmas \ref{cMEq} and \ref{psid} we see
that the extensions $\Theta(E_n)/E$
and $E(\alpha)/E$ satisfy the hypotheses
of Theorem~\ref{conditions}, with $a=tr_n$ and
$m=n-2$.  Therefore there is an automorphism $\Psi$ of
$\Q_p^{alg}$ such that
\begin{equation}
[E(\alpha)\cap\Psi(\Theta(E_n)): E]\ge p^{n-2}.
\end{equation}
Since $E_n=EL_n$ and $\Psi(\Theta(E))=E$, this proves
Theorem~\ref{ext}. \qed

\end{document}